\documentclass[10pt,a4paper]{article}
\usepackage[utf8]{inputenc}
\usepackage{amsmath}
\usepackage{pdfsync}
\usepackage{amsfonts}
\usepackage{amssymb}
\usepackage[left=2cm,right=2cm,top=2cm,bottom=2cm]{geometry}
 \usepackage{epsfig} %For pictures: screened artwork should be set up with an 85 or 100 line screen
\usepackage{graphicx}  \usepackage{epstopdf}%This is to transfer .eps figure to .pdf figure; please compile your paper using PDFLeTex or PDFTeXify.
 \usepackage[colorlinks=true]{hyperref}
\hypersetup{urlcolor=blue, citecolor=red}

\newtheorem{theorem}{Theorem}[section]
\newtheorem{corollary}{Corollary}
\newtheorem{lemma}[theorem]{Lemma}
\newtheorem{proposition}{Proposition}

\newtheorem{definition}[theorem]{Definition}
\newtheorem{remark}{Remark}
   \title{ {Existence and regularity results for some Fully Non Linear singular or degenerate equation }}  
   
\begin{document}
\maketitle
\centerline{\scshape C. O. Ndaw}
\medskip
{\footnotesize
% please put the address of the first author
 \centerline{UMR 8088}
   \centerline{CY Cergy  Paris University  }
   \centerline{ 2 avenue Adolphe Chauvin, Cergy, France    
   }
 % Do not forget to end the {\footnotesize by the sign }

\begin{abstract} 
 In this article we prove existence,  uniqueness and regularity  for the singular equation 
 \begin{eqnarray*}%\label{eqfln}
\begin{cases}
|\nabla u|^{\alpha}(F(D^{2}u)+h(x)\cdot\nabla u)+c(x)|u|^{\alpha}u+p(x)u^{-\gamma}=0\mbox{\ \ in } \  \Omega\\
u>0\ \mbox{\ \ in}\ \Omega ,\mbox{\ \ } u=0\mbox{\ \ \ \ on} \  \partial\Omega
\end{cases}
\end{eqnarray*} when $p$ is some continuous and positive function, $c$ and $h$ are continuous, $\alpha > -1$ and $F$ is Fully non linear elliptic. Some conditions on the first eigenvalue for the  operator $-|\nabla u|^{\alpha}(F(D^{2}u)+h(x)\cdot\nabla u)-c(x)|u|^{\alpha}u$ are required. The results generalize the  well known results of Lazer and McKenna. \smallskip\\
\end{abstract}
} 
\section{Introduction and some useful  tools}\ \\ 

In this work we study the existence, uniqueness and regularity of solution for  the singular equation
\begin{eqnarray}\label{eqfln}
\begin{cases}
|\nabla u|^{\alpha}(F(D^{2}u)+h(x)\cdot\nabla u)+c(x)|u|^{\alpha}u+p(x)u^{-\gamma}=0\mbox{\ \ in } \  \Omega\\
u>0\ \mbox{\ \ in}\ \Omega ,\mbox{\ \ } u=0\mbox{\ \ \ \ on} \  \partial\Omega
\end{cases}
\end{eqnarray}
where $\Omega$ is a bounded  $\mathcal{C}^{2}$ domain, $\alpha>-1$, $p>0$ is  continuous on $\bar{\Omega}$, $c$ and $h$ are continuous on $\bar{\Omega}$. 
$F$ is Fully Non Linear elliptic,  and the  solutions are intented in the viscosity sense, this will be precised below. 

 The results here enclosed generalize the pionneer work of   Lazer and McKenna \cite{LMK}, which consider the case where 
 $h = c = \alpha = 0$ and $F$ is  the Laplacian. In this simple framework, the solutions can be intended in the variational sense, even if  the presence of the singular zero order term $p(x) u^{-\gamma}$ lead the authors in \cite{LMK} to use some tools  usually  employed  in viscosity framework : Existence of  convenient sub- and super-solutions, comparison Theorem, compactness of bounded sequences.  In the present work,  the difficulties are linked both to the singularity/degeneracy  of the  term of derivative of order $1$, and to the singularity of the zero order term.
 
 \medskip

  When $\alpha = 0$  there is a big amount of articles about the existence, the uniqueness and
the regularity of viscosity solutions for the equation $F( D^2 u) = f$,  as the paper of Caffarelli \cite{C}, mainly 
devoted to the ${ \mathcal C}^1$ and higher regularity, the paper of Ishii and Lions \cite{IL}, the user's guide of Crandall, Ishii and Lions \cite{usr}, the 
book of Caffarelli and Cabr\'e \cite{CC},  and  the famous paper  of Ishii \cite{I}. 
   
      The case $\alpha \neq 0$ was introduced by   Birindelli and Demengel in \cite{BD1}, \cite{BD2}, which consider the equations 
    $$ | \nabla u |^\alpha (F( D^2 u) + h(x) \cdot \nabla u )= f,$$ with $f$ and $h$ continuous and bounded.
     They provide a definition of viscosity solution which  fits the case $\alpha <0$  ( note that in that case, the equation is not well defined on a point where  the gradient is zero). In particular this definition can also be   used for the $(\alpha+2)$-Laplacian (when one works with viscosity solutions 
in place of classical  solutions), or to the infinity Laplacian. 
 It  can also  be used for the case $\alpha >0$, note that it is shown in \cite{BDCocv} that the solutions are the same as the classical viscosity solutions in that case. 
 
 More recently Attouchi and Ruosteenoja in \cite{AR}  have refined this definition,  for equations of the form
  $$ |\nabla u|^\gamma (\Delta u+(p-2)\Delta_\infty^N u) = f,$$ with $\gamma>-1$, $p>1$. 
  
Concerning the optimal regularity  expected, the example of 
   $\varphi(r) = r^{2+\alpha \over 1+\alpha}$ which solves 
   $$ | \nabla \varphi |^\alpha \Delta \varphi = \dfrac{(2+\alpha)^{1+\alpha}}{(1+\alpha)^{2+\alpha}}((N-1)\alpha+N)$$
    in the ball $B(0, 1)$, shows that , for $\alpha >0$, this  regularity cannot be better than ${ \mathcal C}^{1,{1 \over 1+ \alpha}}$. 
    
    Coming back to the case $\alpha <0$, a first regularity result is proved in \cite{BDjde}  for solutions of the homogeneous Dirichlet problem. When the operator $F$ is concave or convex, the ${\mathcal C}^2$ regularity holds.
  In   the case $\alpha >0$ and $h = 0$,   Imbert and Silvestre proved a ${\mathcal C}^{1 , \beta}$  \underline{interior} regularity result in \cite{IS}. This result is extended to  a local result "up to the boundary" and to the case where $h\neq 0$, \cite{BDCocv}.
   The interior regularity is precised in \cite{ART}  in the case $\alpha >0$, which can be ${\mathcal C}^{1,{1\over 1+ \alpha}}$ in the cases where  for example $F$ is convex or concave. 
   
   Concerning the equation \eqref{eqfln} in the variational setting,   the work of Crandall Rabinowitz and Tartar  \cite{CRT} extend the results of Lazer McKenna  when the Laplacian is replaced by  some linear uniformly elliptic  operator, positively homogeneous of degree $1$, 
   while,  for the Fully non linear setting,
     in  \cite{FQS} the authors consider  $F(x, u, Du, D^2 u)$  which is homogeneous of degree $1$ with respect to all its arguments and Fully Non Linear Elliptic, with some more general  singularity than $p(x) u^{-\gamma}$. Always when $\alpha = 0$ and $F$ is replaced by some degenerate Pucci's operators ${ \cal P}_k^{\pm}$,  Birindelli and Galise  \cite{BG} proved existence, uniqueness and some regularity result, result extended  to more general   singular zero order term in \cite{V}.

 We now precise the assumptions, and present  the main result. 
We assume that $F$ satisfies the assumptions:

\noindent There exist ellipticity constants $a$ and $A$,  $0<a<A$ so that for any $M$ and $N$ symmetric matrices on $\mathbb{R}^{N}$, $N\geq0$,
\begin{equation}\label{defF}
atr(N)\leq F(M+N)- F(M)\leq Atr(N).
\end{equation} 

We will also assume that $F$ is positively homogenous, i.e 
\begin{equation*}
F(tM)=tF(M)
\end{equation*} 
for all $t>0$.

A well known example of  such operators  are  the Pucci's  operator ${ \cal M}^\pm_{a, A}$ ( ${ \cal P}^+_{ \lambda, \Lambda}$  for  most of the authors),  defined as 

$${ \cal M}^+_{a, A} ( S) = a \sum_{ i, \lambda_i >0} \lambda_i(S) + A \sum_{i, \lambda_i <0} \lambda_i(S),   \ { \cal M}^-_{a, A} (S) = - { \cal M}^+_{a, A} ( -S)$$
 where the $\lambda_i(S) $ are the eigenvalue of the symmetric matrix $S$.

 In the following Theorem $\lambda_1^{\bar c}$ denotes the first eigenvalue for the operator 
$ -| \nabla u |^\alpha (F( D^2 u)+ h(x) \cdot \nabla u)-\bar c(x) |u|^\alpha u$ with the definition precised later. 

 \begin{theorem}\label{thexiF}
Let $\Omega$ be a bounded ${ \cal C}^2$  domain in $\mathbb{R}^{N}$. Suppose that $F$ satisfies \eqref{defF}, is positively homogeneous of degree $1$, that $c$, $h$ , $p$ are continous and bounded, with 
 $p>0$ on $\overline{\Omega}$.
Let us suppose that $\lambda_{1}^{c},\ \lambda_{1}^{\frac{(1+\alpha+\gamma)c}{2+\alpha}}>0,$ then there exists a unique solution to (\ref{eqfln}). In addition $u\ \in\ \mathcal{C}^{1,\beta}(\Omega)$ for some $\beta >0$.
\end{theorem}

  The plan of this paper  is as follows : 
  In  Section 2 we remind the definition of viscosity solutions  adapted to the present context, and recall the  maximum  and comparison principles, and the regularity results needed in the paper.  In Section 3 we prove the existence's result and provide the convenient comparison principle for such equations,   which allows to prove the uniqueness of solutions. In Section 4 we study the regularity of the solutions. 
  
  \section{Background, definitions, and previous existence's and regularity results for singular and degenerate Elliptic equations $| \nabla u |^\alpha F( D^2 u) = f$. } \ \\

We begin to recall the definition of viscosity solutions adapted to the present context. 
We denote by $\mathcal{S}$ the space of symmetric matrices on $\mathbb R^N$. 
Let us define  for $f \in { \mathcal C} ( \Omega)$
\begin{eqnarray}\label{eqgen}
G(x,u,q,X)=|q|^{\alpha}(F(X)+h(x)\cdot q)-f
\end{eqnarray}
where $x\in\mathbb{R}^{N},\ q\in\mathbb{R}^{N},\ u\in\mathbb{R}\ \mbox{and \ }X\in\mathcal{S}$. 
\begin{definition}
A function $u$, upper semicontinuous (USC for short) in $\Omega$ is a viscosity sub-solution for (\ref{eqgen}) (or a solution of $G[u]\geq 0$, a sub-solution of $G[u]= 0$) if whenever $\varphi \in \mathcal{C}^{2}(\Omega)$ and $u-\varphi$ attains a local maximum at $\bar{x}\in \Omega$, then 

1) Either  $\nabla \varphi( \bar x) \neq 0 $  and $$G(\bar{x},u (\bar{x}),D\varphi(\bar{x}),D^{2}\varphi(\bar{x}))\geq 0.$$

2) Or  there exists a ball around $\bar x$ on which $u(x) = u(\bar x)$,  and  
$$-f(\bar x) \geq 0.$$

Similarly, $u$,  lower semicontinous ( LSC for short) in $\Omega$ is a viscosity super-solution for  (\ref{eqgen}) (or a solution of $G[u]\leq 0$, a super-solution of $G[u]=0$) if whenever $\varphi \in \mathcal{C}^{2}(\Omega)$ and $u-\varphi$ attains a local minimum at $\bar{x} \in \Omega$, and $\nabla \varphi(\bar x)\neq 0$
 then $$G(\bar{x},u (\bar{x}),D\varphi(\bar{x}),D^{2}\varphi(\bar{x}))\leq 0.$$
  If $u$ is locally constant around $\bar{x}$
  $$-f( \bar x) \leq 0.$$
  
Of course, $u \in \mathcal{C}(\Omega)$ is a viscosity solution of (\ref{eqgen}) (or a solution of $G[u]= 0$) if $u$ is both a viscosity sub-solution and a viscosity super-solution.
\end{definition}
Let $u\in \mathcal{C}(\Omega)$. We define the superjet $J^{2,+}u(x)$ and the subjet $J^{2,-}u(x)$ of the second order.
\begin{definition}
\begin{eqnarray*}
J^{2,+}u(x)=\left\lbrace (p,X)\in \mathbb{R}^{N}\times\mathbb{S}\right.,  {\rm\ so\ that} & u(x)\leq u(\bar{x}) +<p,x-\bar{x}> +\frac{1}{2}<X(x-\bar{x}),x-\bar{x}>\\
&\left. +o(| x-\bar{x}|^{2})\right\rbrace.\\
J^{2,-}u(x)=\left\lbrace (p,X) \in \mathbb{R}^{N}\times\mathbb{S}\right.{\rm,\ so\ that}  & u(x)\geq u(\bar{x}) +<p,x-\bar{x}>+\frac{1}{2}<X(x-\bar{x}),x-\bar{x}>\\
&\left. +o(|x-\bar{x}|^{2})\right\rbrace.
\end{eqnarray*}
\end{definition} 
 More useful are the closed superjet and closed subjet,
 $$\overline{J^{2,+}u(\bar x)} =  \{(\bar p, \bar  X),\ \exists x_n, \ x_n \rightarrow \bar x, \ {\rm and} \  (p_n , X_n) \in J^{2,+}u(x_n), \ (p_n, X_n) \rightarrow ( p, X)\}.$$
$\overline{ J^{2,-} u( \bar x)}$ being defined in an obvious symmetric manner. 

In the definition of viscosity solutions the test functions can be substituted by the elements of the semi-jets in the sense that in the definition above one can restrict to  the functions $\varphi$ defined by $$\varphi(x)=u(\bar{x})+<p,x-\bar{x}>+\frac{1}{2}<X(x-\bar{x}),x-\bar{x}>$$ with $(p,X) \in J^{2,-}u(\bar{x})$ when $u$ is a super-solution and $(p,X) \in J^{2,+}u(\bar{x})$ when u is a sub-solution.
\ \\
A key tool for  the existence's results is  the comparison principle, \cite{BD1}.

\begin{theorem}\label{Comparisontheorem}
Let $\Omega$ be a bounded ${\mathcal C}^2$ domain in $\mathbb R^N$. Let $k$  be  continuous with respect to its variables.   
Suppose that $u$ is a USC solution of 
\begin{eqnarray*}
|\nabla u|^{\alpha}(F(D^{2}u)+h(x)\cdot\nabla u)-k(x, u)\geq f
\end{eqnarray*}
and $v$ a LSC solution of 
\begin{eqnarray*}
|\nabla v|^{\alpha}(F(D^{2}v)+h(x)\cdot\nabla v)-k(x, v)\leq g
\end{eqnarray*}
 with $f,g \in \mathcal{C}(\Omega)$, $f\geq g$ and $u\mapsto k(x, u)$  increasing, or $f>g$ and $u\mapsto k(x,u)$ non decreasing.\\
 Then if $u\leq v$ on $\partial \Omega$, $u\leq v$ in $\Omega$.
\end{theorem}

This comparison theorem, and the  construction of sub- and super-solutions which are zero on the boundary,  together with  a mere adaptation to our context of Perron's method, permit to prove the existence of solutions  of the Dirichlet problem, for 
$$ \left\{ \begin{array}{lc} 
| \nabla u |^\alpha ( F( D^2 u)+ h(x) \cdot \nabla u) -k (x,  u) = g& {\rm in}\ \Omega \\
 u = 0 & {\rm on} \ \partial \Omega. \end{array}\right.$$
  Since we are dealing with positive solutions, we will many times use the strong maximum principle :

\begin{theorem}\label{Pmaxf} Let $\Omega$ be a bounded ${ \cal C}^2$  domain in $\mathbb R^N$. Suppose that $h$ and $c$ are continuous and let $u\geq0$ satisfies
\begin{eqnarray*}
|\nabla u|^{\alpha}(F(D^{2}u)+h(x) \cdot \nabla u) + c(x)u^{1+\alpha}\leq 0\mbox{\ \ \ \ \ \ } in\ \Omega. 
\end{eqnarray*}
Then one has
\begin{eqnarray*}
u>0 \mbox{\ \ \ \ \ or\ \ \ \ \ }u\equiv 0.
\end{eqnarray*}
\end{theorem}

 For some of the results enclosed we will need  "Hopf boundary principle", say the fact that near the boundary, the gradient of positive solutions cannot be zero in a neighborhood of the boundary : 
\begin{theorem}\label{PHopf} Let $\Omega$ be a bounded ${ \cal C}^2$  domain in $\mathbb R^N$. Let $x_{0} \in \partial\Omega$ and $\partial\Omega$ satisfy the  interior sphere condition at $x_o$. Let $\overrightarrow{n}$ denotes the inner normal to $\partial\Omega$ at $x_{0}$. If $u>0$ in $\Omega $  and 
\begin{eqnarray*}
|\nabla u|^{\alpha}F(D^{2}u)+h(x)\cdot\nabla u|\nabla u|^{\alpha}+c(x)u^{1+\alpha}\leq 0\mbox{\ \ \ \ \ \ } in\ \Omega \mbox{\ \ and\ \ } u(x_{0})=0  
\end{eqnarray*}
then
\begin{eqnarray*}
\dfrac{\partial u}{\partial \overrightarrow{n}}(x_{0})>0.
\end{eqnarray*}
in the sense
\begin{eqnarray*}
\underline{\lim}_{h\rightarrow 0, h>0}{\dfrac{u(x_{0}+h\overrightarrow{n})-u(x_{0})}{h}>0}.
\end{eqnarray*}
\end{theorem}
\begin{remark} As an easy consequence, if $u$ is a solution of the equation  which is $\mathcal{C}^{1}$ up to the boundary,  we have for some $\kappa >0$ , for any $x_o\in \partial \Omega$ such that $u(x_{0})=0$,
\begin{eqnarray*}
\vert\dfrac{\partial u}{\partial \overrightarrow{n}}(x_{0})\vert \geq \kappa\ \mbox{\rm  \ \ \ and then\ \ \ \ \ }\vert \nabla u\vert \geq \kappa  \mbox{\rm  \ \ \ in a neighborhood of \ \  }x_o.
\end{eqnarray*}
\end{remark}

We now recall the Lipschitz estimates between sub-and super-solutions (See \cite{IL} for example in the case $\alpha = 0$).

 \begin{theorem}\label{lip} 
1) Let $u$ be USC such that 
\begin{eqnarray*}
|\nabla u|^{\alpha}(F(D^{2}u)+ h(x) \cdot \nabla u) \geq  f\ \mbox{\ \ \ \ \ } {\rm in}\ B(0,1)
\end{eqnarray*}
and $v$ is LSC and satisfies 
\begin{eqnarray*}
|\nabla v|^{\alpha}(F(D^{2}v)+ h(x) \cdot \nabla v) \leq  g \ \mbox{\ \ \ \ \ }{\rm  in} \ B(0,1)
\end{eqnarray*}
with $f$ and $g$ continuous and bounded. 

Then $\forall\ r\in (0,1)$, there exists $\ L_{r}>0$ such that $\forall\ (x, y) \in ({B}(0,r))^{2}$,
\begin{eqnarray*}
u(x)-v(y)\leq \sup ( u-v) +  L_{r}|x-y|.
\end{eqnarray*}
Here $L_r$ depends on $r,$ $\sup u-\inf v$, $|f|_\infty$ and $|g|_\infty $.

 In particular any  solution  of $|\nabla v|^{\alpha}(F(D^{2}v)+ h(x) \cdot \nabla v) =  f $, with $f$ continuous and bounded,  is  locally Lipschitz continuous. 
 
 2) Suppose that $u$ is a solution of the Dirichlet problem 
 
$$\left\{ \begin{array}{lc}
|\nabla u|^{\alpha}(F(D^{2}u)+ h(x) \cdot \nabla u) =   f, \mbox{\ \ \ \ \ } &f \ {\rm in} \ B(0,1),\\
u=0 & {\rm on} \ \partial (B(0,1)).
\end{array}\right. 
$$
 Then  $u$  is Lipschitz continuous, with some Lipschitz constant depending on $|f|_\infty$ and $|u|_\infty$. 
 
\end{theorem}

As a corollary  one easily has by using the definition of viscosity sub- and super-solutions

   \begin{theorem}\label{Coseqlips}
   1) Let $\Omega$ be a bounded $\mathcal{C}^{2}$ domain in $\mathbb{R}^{N}$. Let $u_n$  be a viscosity solution of 
    $$| \nabla  u_n|^\alpha (F(  D^2 u_n) + h(x) \cdot \nabla u_n)= f_n$$
    in $\Omega$
 and suppose that $f_n$ converges locally  uniformly to $f$,  and $u_n$ is locally uniformly bounded, then one can extract from $(u_n)_n$ a subsequence such that this subsequence converges locally uniformly towards a solution of 
   $$ | \nabla  u|^\alpha (F( D^2 u)+ h(x) \cdot \nabla u) = f.$$
  2) If moreover $u_n= 0$ on the boundary and $f_{n}$ converges  uniformly to $f$, $u_n$ converges uniformly up to a subsequence towards $u$ on $\overline{ \Omega}$. 
\end{theorem}

\medskip

 Let us now  make precise the   ${\mathcal C}^{1, \beta}$ regularity results   for the solutions of 
 $$| \nabla u |^\alpha (F( D^2 u) + h(x) \cdot \nabla u) = f,$$ when $f$ is continuous in $\Omega$. 
  The first result obtained in \cite{BDjde} is 
  
  \begin{theorem}
Let $\Omega$ be a bounded ${ \cal C}^2$  domain in $\mathbb R^N$.  Let $\alpha \in ]-1,0]$.  Let $f$ and $h$ be continuous on $\overline{\Omega}$. There exists $\beta$  so that for any $u$ solution of 
    $$\left\{\begin{array}{lc}
     | \nabla u |^\alpha (F( D^2 u) + h(x) \cdot \nabla u) = f& {\rm in} \ \Omega\\
      u = 0 & {\rm on} \ \partial \Omega, \end{array}\right. $$
     $u \in {\mathcal C}^{1, \beta} ( \overline{ \Omega})$.
      \end{theorem}
      
       The case $\alpha >0$ is treated  in \cite{IS}, for $h = 0$, the precise result is an interior estimate: 
       
       \begin{theorem}
        Suppose that $\alpha \geq 0$. There exists $\beta \in ]0,1[$, such that for all $r \in ]0,1[$, there exists $C_r$ so that for any $f\in {\mathcal C} ( B(0,1))$, and for any $u$ solution of 
        $$ | \nabla u |^\alpha F( D^2 u) = f$$ in $B(0,1)$, one has that
        $u$ is $ {\mathcal C}^{1, \beta} ( B(0,r))$, with 
        $$|u|_{{\mathcal C}^{1, \beta} ( B(0, r))} \leq C_r ( |u|_\infty + |f|_\infty^{1\over 1+ \alpha}).$$
        \end{theorem}
        
        To complete this interior estimate, a  ${\mathcal C}^{1, \beta} ( \overline{ \Omega})$ regularity result  is obtained in  \cite{BDCocv}  : 
        \begin{theorem}
        
        Suppose that $\Omega$
 is a bounded ${\mathcal C}^2$ domain in $\mathbb R^N$, that $f$ and $h$ are continuous on $\overline{ \Omega}$. 
 Let $\varphi \in {\mathcal C}^{1, \beta_o} ( \partial \Omega)$. There exists $\beta \leq \inf (  \beta_o, {1\over 1+ \alpha})$, and some constant $C$  so that  for any $u $ solution of 
  $$\left\{\begin{array}{lc}
     | \nabla u |^\alpha (F( D^2 u) + h(x) \cdot \nabla u) = f& {\rm in} \ \Omega\\
      u = \varphi  & {\rm on} \ \partial \Omega, \end{array}\right. $$
      one has that $u\in {\mathcal C}^{1, \beta}$ on $\overline{\Omega}$  and 
       $$|u|_{{\mathcal C}^{1, \beta}  ( \overline{ \Omega})}\leq C ( |\varphi|_{ {\mathcal C}^{1, \beta_o}( \partial \Omega)} + |f|_\infty^{1\over 1+\alpha}).$$
       \end{theorem}
        
         Other authors  provide  precises bounds on   $\beta$, the interested reader can see     \cite{ART}.

 We now precise the definition of the first  demi eigenvalues and eigenfunctions,  \cite{BD2}.   They are defined on the model of \cite{BNV} : 
  
   $$\lambda_1 ^{+, c} = \sup \{ \lambda, \exists \varphi >0, |\nabla \varphi |^\alpha ( F( D^2 \varphi) +  h(x) \cdot \nabla \varphi) + ( c(x) + \lambda) |\varphi |^\alpha \varphi \leq 0\}. $$
   
      It is clear that $\lambda_1^{+,c} $ exists and $\lambda_1^{+,c}\geq -|c|_\infty$. 
      
    Some precise estimates depending on the larger ball contained in $\Omega $,  and the smallest ball  containing it,  can be found in \cite{BD2}. 
    Even if this article will only need  $\lambda_1^{+, c}$,  we give the definition of the other demi eigenvalue, say:       $$\lambda_1 ^{-, c} = \sup \{ \lambda, \exists \varphi <0, 
      |\nabla \varphi |^\alpha ( F( D^2 \varphi) +  h(x) \cdot \nabla \varphi) + ( c(x) + \lambda) |\varphi |^\alpha \varphi \geq 0\}.$$
    Note that in the case where $F$ is odd, the two eigenvalues coincide.

     We have the following maximum and minimum principle "under the first eigenvalues"
     
      \begin{theorem}Let $\Omega$ be a bounded ${ \cal C}^2$  domain in $\mathbb R^N$. Under the previous assumptions on $F$, $h, c$, $\alpha$,  
       suppose that $\tau < \lambda_1^{+, c}$ and that $u$, USC is a sub-solution of
       $$ | \nabla u |^\alpha(  F( D^2 u) + h(x) \cdot \nabla u) +( c(x) + \tau) |u|^\alpha u  \geq 0 \ {\rm in} \ \Omega$$
        and $u\leq 0$ on $\partial \Omega$. Then 
        $u\leq 0$ in $\Omega$.
        
          Suppose that $\tau < \lambda_1^{-, c}$ and that $u$, LSC is a super-solution of
       $$ | \nabla u |^\alpha(  F( D^2 u) + h(x) \cdot \nabla u) +( c(x) + \tau) |u|^\alpha u  \leq 0 \ {\rm in} \ \Omega$$
        and $u\geq 0$ on $\partial \Omega$. Then 
        $u\geq 0$ in $\Omega$.
        \end{theorem}
        
         This Theorem allows to prove the existence of a positive eigenfunction for $\lambda_1^{+, c}$ and a negative one for $\lambda_1^{-, c} $: 
         
          \begin{theorem}Let $\Omega$ be a bounded ${ \cal C}^2$  domain in $\mathbb R^N$.Under the previous assumptions on $F$, $h, c$, $\alpha$,  
          there exists a positive eigenfunction associated to $\lambda_1^{+, c}$, more precisely 
          $$ \left\{\begin{array}{lc}
           | \nabla \varphi_1^+|^\alpha (F( D^2 \varphi_1^+)+ h(x) \cdot \nabla \varphi_1^+) 
           + ( c(x) + \lambda_1^{+, c} )(\varphi_1^+)^{1+ \alpha} = 0  &{\rm in} \ \Omega\\
            \varphi_1^+ = 0 & {\rm on } \ \partial \Omega. 
            \end{array}\right.$$
             There exists a negative  eigenfunction associated to $\lambda_1^{+, c}$, more precisely 
          $$ \left\{\begin{array}{lc}
           | \nabla \varphi_1^-|^\alpha (F( D^2 \varphi_1^-)+ h(x) \cdot \nabla \varphi_1^-) + ( c(x) + \lambda_1^{-, c}) |\varphi_1^-|^{ \alpha}\varphi_1^- = 0& {\rm in} \ \Omega\\
            \varphi_1^- = 0 & {\rm on } \ \partial \Omega.
            \end{array}\right.$$
            \end{theorem}
            In the following we will drop the exponent $+$ since we will only use $\lambda_1^{+, c}$. 
 
\section{Existence and uniqueness of  viscosity solutions}
\subsection{Existence of  viscosity sub- and super-solutions : Proof of Theorem \ref{thexiF}}\ \\

Let  $\beta_{c,\alpha,\gamma}=\dfrac{(1+\alpha+\gamma)c}{2+\alpha}$.  
 We begin with some remark about the validity of assumption $\lambda_1^c >0$, assumed in the existence's Theorem: 
\begin{remark}
If we take for example $|c|_{\infty}<\lambda_{1}^0$, which is known to be $>0$ one has $\lambda_{1}^{c}>0$. The same is true for $|\beta_{c,\alpha,\gamma}|_{\infty}<\lambda_{1}^{0}$. \
\end{remark}
We begin to exhibit a sub- and a super-solution.  Let $\phi_{1}$ be an eigenfunction for $\lambda_{1}^{\beta_{c,\alpha,\gamma}}$. We  first treat the case $\gamma >1$ : 
\begin{proposition}
Let assume $\gamma>1$. Let $t=\dfrac{2+\alpha}{1+\alpha+\gamma}$. There exist $b_{i},\ i=1,\ 2$ so that $\psi_{i}=b_{i}\phi_{1}^{t}$ are respectively sub- and super-solutions.
\end{proposition}
\textbf{Proof.}
 We do the sub-solution case. Let $\psi=b\phi_{1}^{t}$. Then $\nabla  \psi = bt \phi_1^{t-1} \nabla \phi_1$. 
\begin{eqnarray*}
D^{2}\psi &=&bt(t-1)\phi_{1}^{t-2}\nabla\phi_{1}\otimes\nabla\phi_{1}+bt\phi_{1}^{t-1}D^{2}\phi_{1}.  
\end{eqnarray*}
Then 

\begin{eqnarray*}
|\nabla \psi|^{\alpha}(F(D^{2}\psi) &+&h(x)\cdot\nabla \psi)+c(x)\psi^{1+\alpha}\\  
&\geq& Ab^{1+\alpha}t^{1+\alpha}(t-1)\phi_{1}^{t-2+(t-1)\alpha}|\nabla\phi_{1}|^{2+\alpha }\\ 
&+&b^{1+\alpha}t^{1+\alpha} \phi_{1}^{(t-1)(1+\alpha)}\left(|\nabla \phi_{1}|^{\alpha}(F(D^{2}\phi_{1})+h(x)\cdot\nabla\phi_{1})+\dfrac{c}{t^{1+\alpha}}\phi_{1}^{1+\alpha}\right)\\ 
&=&  (bt)^{1+\alpha} \phi_1^{(t-1)(1+\alpha)} \left( | \nabla \phi_1|^\alpha (F( D^2 \phi_1) + h(x) \cdot \nabla \phi_1)+( { c \over t^{1+\alpha}} +\lambda_1^{\beta_{c,\alpha,\gamma}} )\phi_1^{1+\alpha} \right) \\ 
  &-& b^{1+\alpha}  t^{1+\alpha}\phi_1^{(t-1) \alpha + t-2}( A  (1-t)| \nabla \phi_1|^{2+\alpha} +  \lambda_1^{\beta_{c,\alpha,\gamma}} \phi _1^{2+\alpha} )\\ 
&= &-b^{1+\alpha}t^{1+\alpha}\phi_{1}^{t-2+(t-1)\alpha}\left( (1-t)A|\nabla\phi_{1}|^{2+\alpha}+\lambda_{1}^{\beta_{c,\alpha,\gamma}}\phi_{1}^{2+\alpha}\right) 
\end{eqnarray*}
and then if we denote
\begin{eqnarray*}
q^{A}(x, b)=b^{1+\alpha}t^{1+\alpha}\phi_{1}^{t-2+(t-1)\alpha}\left( (1-t)A|\nabla\phi_{1}|^{2+\alpha}+\lambda_{1}^{\beta_{c,\alpha,\gamma}}\phi_{1}^{2+\alpha}\right),
\end{eqnarray*}
one has
\begin{eqnarray*}
|\nabla \psi|^{\alpha}(F(D^{2}\psi)+h(x)\cdot\nabla \psi)+c(x)\psi^{1+\alpha}+q^A(x,b)\geq 0.
\end{eqnarray*}
 Note that using analogous computations, one has 
 \begin{eqnarray*}
|\nabla \psi|^{\alpha}(F(D^{2}\psi)+h(x)\cdot\nabla \psi)+c(x)\psi^{1+\alpha}+q^a(x,b)\leq 0
\end{eqnarray*}
where $q^a $ is defined by replacing $A$ by $a$ in the definition of $q^A$. 

\textbf{Claim:} one has the existence of positive constants $d_{i},\ i=1,2$ so that
\begin{eqnarray}\label{ineqd1d2}
d_{2}\leq A(1-t)|\nabla\phi_{1}|^{2+\alpha}+\lambda_{1}^{\beta_{c,\alpha,\gamma}}\phi_{1}^{2+\alpha} \leq d_{1}.
\end{eqnarray}
 
Let us admit for a while the claim, and   let us take 
\begin{eqnarray*}
b_{1}=\left( \dfrac{\min{p}}{d_{1}t^{1+\alpha}}\right)^{\frac{1}{1+\alpha+\gamma}},
\end{eqnarray*}
then $q^{A}(x,b_{1})\leq p(x)\psi_{1}^{-\gamma}$  and 

\begin{eqnarray*}
|\nabla (b_{1}\phi)^{t})|^{\alpha}(F(D^{2}(b_{1}\phi^{t}) )&+&h(x)\cdot\nabla (b_{1}\phi^{t}))\\
&+&c(x)(b_{1}\phi_{1}^{t})^{1+\alpha}+p(x)(b_{1}\phi_{1}^t)^{-\gamma}\\
&\geq &|\nabla (b_{1}\phi_{1})^{t})|^{\alpha}(F(D^{2}(b_{1}\phi_{1})^{t})+h(x)\cdot\nabla (b_{1}\phi_{1})^{t}))\\
 & &+c(x)(b_{1}\phi_{1}^{t})^{1+\alpha}+q^{A}(x,b_{1}).
\end{eqnarray*}
Then with that choice of $b_1$,  $\psi_{1}=b_{1}\phi_{1}^{t}$ is a sub-solution of (\ref{eqfln}).\\
In  the same  manner using the left hand  side inequality of (\ref{ineqd1d2}) and taking $$
b_{2}=\left( \dfrac{\max{p}}{d_{2}t^{1+\alpha}}\right)^{\frac{1}{1+\alpha+\gamma}},
$$
  $\psi_{2}=b_{2}\phi_{1}^{t}$ is a super-solution of (\ref{eqfln}). Note that  $b_{1}<b_{2}$ and then $\psi_1 \leq \psi_2$.
  
We now prove  claim \eqref{ineqd1d2}:  
Since $\phi_{1}$ is in $\mathcal{C}^{1}(\overline{\Omega})$ the inequality on the right side  of \eqref{ineqd1d2} is obvious.

To prove the left hand side inequality, 
near the boundary, using Hopf principle and  the fact that  $\phi_{1}$ is of class $\mathcal{C}^{1}$, $\exists\ \delta,\ d(x,\partial\Omega)<\delta \ $ $ \Rightarrow|\ \nabla \phi_{1}|>m$ (anywhere else one will use $|\nabla \phi_{1}|\geq  0$).\\
Now, using  $\phi_{1}>0$,  
  $\exists\ \tilde{m}$ such that $\phi_{1}\geq\tilde{m}$ on $d(x,\partial\Omega)\geq \delta$. From all this we derive that 
\begin{align*}
c_{1}|\nabla\phi_{1}|^{2+\alpha}+c_{2}\phi_{1}^{2+\alpha}\geq \min{(c_{1},c_{2})}\min{(m^{2+\alpha},\tilde{m}^{2+\alpha})}.
\end{align*}
Taking
\begin{align*}
d_{2}= \min{(c_{1},c_{2})}\min{(m^{2+\alpha},\tilde{m}^{2+\alpha})}
\end{align*}
we have the left-side inequality of \eqref{ineqd1d2}.
$\hfill\square$
\ \\

We now treat the case $\gamma<1$.  We begin to introduce  a regularized problem depending on  some  parameter $\delta >0$,
\begin{equation}\label{eqreg}
\left\lbrace\begin{array}{lc}|\nabla u|^{\alpha}(F(D^{2}u)+h(x)\cdot\nabla u)+c(x)u^{1+\alpha}+p(x)(u+\delta)^{-\gamma}=0 & {\rm in} \ \Omega, \\
u=0 & {\rm on } \ \partial \Omega .
\end{array}
\right.
\end{equation}
\begin{lemma}
 Let $\gamma>0$,  and $c$ be so that  $\lambda_{1}^{c}>0$, let $\psi_{1}$ be some  positive eigenfunction for $\lambda_{1}^{c}$. Then, there exist $\varepsilon_{0}$ and $\delta_{0}$ so that for $\varepsilon<\varepsilon_{0}$ and $\delta\in [0, \delta_{0}],$   then $u_{\star}=\varepsilon\psi_{1}$ is a sub-solution of\eqref{eqreg}.

\end{lemma}
\textbf{Proof.}

Let $u_{\star} = \varepsilon \psi_1$. Take
\begin{eqnarray*}
\delta_{0}= \dfrac{1}{2^{\frac{\gamma}{1+\alpha+\gamma}}}\left( \frac{\min{p}}{ \lambda_{1}^{c}}\right)^{\dfrac{1}{1+\alpha+\gamma}} \mbox{\ \ and\ \ }\varepsilon_{0}= \dfrac{1}{2^{\frac{\gamma}{1+\alpha+\gamma}}|\psi_1|_{\infty}}\left( \frac{\min{p}}{\lambda_{1}^{c}}\right)^{\dfrac{1}{1+\alpha+\gamma}}.
\end{eqnarray*}
Then,  by an easy computation,  one has for all $\varepsilon<\varepsilon_{0}$ and $\delta<\delta_{0}$, 
$$| \nabla u_{\star}|^{\alpha}(F(D^{2}u_{\star})+h(x)\cdot\nabla u_{\star})+c(x)u_{\star}^{1+\alpha}+p(x)(u_{\star}+\delta)^{-\gamma}\geq 0. 
$$
$\hfill\square$
\ \\

\begin{proposition} \label{existsub}
Suppose that  $\gamma<1$,  and $s<1$ sufficiently close to $1$ in order 
 that $\lambda_{1}^{cs^{-(1+\alpha)}}>0 $. \\
If $\psi_{2}$ is an eigenfunction corresponding to $\lambda_{1}^{cs^{-(1+\alpha)}}$, then there exists $d$ great enough  in order that $u^{\star}=d\psi_{2}^{s}$ is a super-solution of (\ref{eqfln}).
\end{proposition}
\begin{remark} The fact that $\lambda_{1}^{cs^{-(1+\alpha)}}>0 $ when $s$ is sufficiently close to $1$ is justified by the continuity result in the Proposition \ref{prop35} below.
\end{remark}
\ \\
\textbf{Proof.}  of Proposition \ref{existsub}
Let $u^{\star} =  d\psi_{2}^{s}$. We have
\begin{eqnarray*}
|\nabla u^{\star}|^{\alpha}(F(D^{2}u^{\star})&+&h(x)\cdot\nabla u^{\star})+ c(u^{\star})^{1+\alpha}+p(x)(u^{\star})^{-\gamma}\\
&\leq& d^{1+\alpha}s^{1+\alpha}\psi_{2}^{s-2+(s-1)\alpha}\left( (s-1)|\nabla\psi_{2}|^{2+\alpha}a\right.\\ 
&&\left.+|\nabla \psi_{2}|^{\alpha}(F(D^{2}\psi_{2})+h(x)\cdot\nabla \psi_{2})\right)+c(x)(d\psi_{2}^{s})^{1+\alpha}+p(x)d^{-\gamma}\psi_{2}^{-\gamma s}\\
 &=&  (ds)^{1+\alpha} \psi_2^{(s-1)(1+\alpha)} \left( |\nabla \psi_2|^\alpha F( D^2 \psi_2) + h(x) \cdot \nabla \psi_2+( { c \over s^{1+\alpha}} +\lambda_1^{c s^{-(1+\alpha)}})\psi_2^{1+\alpha} \right) \\
  &-& d^{1+\alpha}  s^{1+\alpha}  \psi_2^{(s-1) \alpha + t-2}( a  (1-s)| \nabla \psi_2|^{2+\alpha} + \lambda_1^{c s^{-(1+\alpha)}}\psi _2^{2+\alpha} )+p(x)d^{-\gamma}\psi_{2}^{-\gamma s}\\ 
&= &-d^{1+\alpha}s^{1+\alpha}\psi_{2}^{s-2+(s-1)\alpha}\left( (1-s)|\nabla\psi_{2}|^{2+\alpha}a+\lambda_1^{\frac{c}{s^{1+\alpha}}}\psi_{2}^{2+\alpha}\right)+p(x)d^{-\gamma}\psi_{2}^{-\gamma s}.\\
\end{eqnarray*}
Since $\gamma<1$ one has $\dfrac{2+\alpha}{1+\alpha+\gamma}>1$, and then if $s<1$, one has $-s(\gamma+1)+ (1-s)\alpha+2>0$.\\
Then denoting
$$\kappa = \min_{x\in \overline{ \Omega}}\left(  (1-s)|\nabla\psi_{2}|^{2+\alpha}a+\lambda_1^{\frac{c}{s^{1+\alpha}}}\psi_2^{2+\alpha}\right)$$
(the existence of $\kappa$ can be proved in the same manner as the existence of $d_2$ in \eqref{ineqd1d2})  and 
assuming $d$ large enough in order that 
\begin{eqnarray*}
d\geq\left( \dfrac{|p|_{\infty}(\max \psi_{2}(x))^{-s(\gamma+1)+ (1-s)\alpha+2}}{\kappa s^{1+\alpha}}\right) ^{1\over {\alpha+\gamma+1}}.
\end{eqnarray*}
We have
\begin{eqnarray*}
|\nabla u^{\star}|^{\alpha}(F(D^{2}u^{\star})+h(x)\cdot\nabla u^{\star})+ c(u^{\star})^{1+\alpha}+p(x)(u^{\star})^{-\gamma}\leq 0.
\end{eqnarray*}
So $u^{\star}$ is a super-solution of  equation  \eqref{eqfln}.
$\hfill\square$
\ \\
Furthermore for $\varepsilon$ chosen small and and $d$ large  we have that $u_\star  < u^\star$. 
Indeed let $\psi_1$ be a positive eigenfunction  for $\lambda_1^c$ and $\psi_2$ a positive eigenfunction for $\lambda_1^{c s^{-(1+\alpha)}}$, by the results in \cite{BD1}, there exist  positive constants 
      $c_i$ $i = 1, \cdots, 4$ so that  
      $$ c_1d(x, \partial \Omega) \leq \psi_1 \leq c_2 d(x, \partial \Omega) ,\  c_3 d(x, \partial \Omega) \leq \psi_2 \leq c_4 d(x, \partial \Omega) $$
       and then  taking $\varepsilon$  small enough and $d$  large enough so that 
        $$ \varepsilon c_2 ({\rm diam} \ \Omega)^{1-s} < d c_3^s$$
         one gets that 
         $$\varepsilon  \psi_1 \leq \varepsilon c_2 d(x, \partial \Omega) <  d c_3^s d(x, \partial \Omega)^s \leq d \psi_2^s. $$

\begin{proposition}\label{prop35} The map $s \mapsto \lambda_1^{cs^{-(1+\alpha)}}$ is continuous. More generally,  if $c_{n}$ converge uniformly to $c$ we have  $\lim \lambda_1^{c_n}  = \lambda_1^c$.
\end{proposition} 
\textbf{Proof.} 
On one hand we have  $\limsup \lambda_{1}^{c_n} \leq \lambda_{1}^c$.
 Indeed, let $\lambda = \limsup  \lambda_{1}^{c_n}$, there exists $\varphi_n$ so that 
   $\varphi_n >0$, $| \varphi_n|_\infty = 1$  and $$| \nabla \varphi_n|^\alpha (F( D^2 \varphi_n)+h(x)\cdot\nabla\varphi_n) + ( c_n+ \lambda_1^{c_n} ) \varphi_n =  0.$$
   Using the uniform Lipschitz estimate in   Theorem \ref{lip},  $(\varphi_n)$ is uniformly Lipschitz continuous, then one can extract from $( \varphi_n)_n$ a subsequence which converges uniformly on $\overline{ \Omega}$ toward some function $\varphi$, in particular $(c_n+ \lambda_1^{c_n}) \varphi_n$ converges uniformly toward $(c+\lambda) \varphi$, and using  Theorem  \ref{Coseqlips}, one gets that $\varphi$ satisfies 
   $$| \nabla \varphi|^\alpha (F( D^2 \varphi)+h(x)\cdot\nabla\varphi) + ( c+ \lambda) \varphi  =  0.$$
   By the strong maximum principle, since $| \varphi|_\infty = 1$,  $\varphi>0$ in $\Omega$ and then
 by the definition of $\lambda_1^c $, 
   $$ \lambda_{1}^c \geq  \lambda.$$     On the other hand let $\lambda < \lambda_1^c$, by the existence's result in \cite {BD1} there exists $\psi$ so that $\psi >0$ and 
    $$| \nabla \psi |^\alpha (F( D^2 \psi) +h(x)\cdot\nabla\psi)+ (\lambda+ c) \psi = -1.$$
    Let then $N$ so that for $n > N$ 
    $ |c_n-c|| \psi|_\infty < {1\over 2}$. Then 
     for such $n$
     $$| \nabla \psi |^\alpha( F( D^2 \psi) +h(x)\cdot\nabla\psi)+ (\lambda+ c_n) \psi \leq  {-1\over 2}$$ and then  for this range of values of $n$, 
      $\lambda_{1}^{c_n} \geq \lambda$. Since $\lambda$ is arbitrary less than $\lambda_1^c$ one gets 
      $ \liminf \lambda_1^{c_n} \geq \lambda_1^c$.
$\hfill\square$
\ \\

\begin{remark}
The function $u^{\star}$  is also a super-solution  of the regularized  problem (\ref{eqreg}). 
\end{remark}

\subsection{Proof of the existence's Theorem}\ \\

Let us now prove  the existence result in Theorem \ref{thexiF}. We will proceed  in two steps.
\subsubsection*{Step1:}
Let $u_{\star}$ and $u^\star$ be respectively the  sub- and super-solutions  (both for the regularized problem) constructed above, with $u_{\star}\leq u^\star$.
Let $k$ be defined by 
\begin{eqnarray*}
k>\max\left\lbrace\dfrac{\gamma}{1+\alpha}\dfrac{|p|_\infty}{\delta^{\alpha+\gamma+1}};|c|_\infty\right\rbrace. 
\end{eqnarray*}
Then
the functions 
\begin{eqnarray*}
f_{\delta}(x,u)=-k(u+\delta)^{1+\alpha}-p(x)(u+\delta)^{-\gamma}
\end{eqnarray*}
and
\begin{eqnarray*}
g_{\delta}(x,u)=c(x)(u+  \delta( \alpha))^{1+\alpha}-k(u+\delta)^{1+\alpha}
\end{eqnarray*}
where $$\delta ( \alpha) = \left\{ \begin{array}{c}\delta \ {\rm if} \ \alpha <0\\
 0\  {\rm   if\ not}, 
 \end{array} \right.$$
are  decreasing  with respect to  $u>0$. 

\noindent In the sequel will suppose $\alpha \geq 0$, the case $\alpha <0$ is left to the reader.

\noindent Let us consider the sequence  $\{w_{n}\}$, defined in a recursive way by
\begin{equation*}
\left\{\begin{array}{lc}
|\nabla w_{n}|^{\alpha}(F(D^{2}w_{n})+h(x)\cdot\nabla w_{n})+c(x)|w_{n}|^{\alpha}w_{n}-k(w_{n}+\delta)^{1+\alpha}=f_{\delta}(x,w_{n-1})&{\rm  in} \ \Omega\\
w_{n}=0\  \mbox{on  }\partial\Omega&
\end{array}\right.,
\end{equation*}
with $w_{0}=u_{\star}$.\\
We will prove that for all $n\in \mathbb{N}$, $u_{\star}\leq w_{n}\leq w_{n+1}\leq u^{\star}$. Let us show  by induction that $\{w_{n}\}$ is non decreasing.\\
To prove that $w_{1}\geq w_{0}$ note that :\\

\begin{eqnarray*}
|\nabla w_{0}|^{\alpha}(F(D^{2}w_{0})+h(x)\cdot\nabla w_{0})+c(x)|w_{0}|^{\alpha}w_{0}-k(w_{0}+\delta)^{1+\alpha}\geq  f_{\delta}(x,w_{0})
\end{eqnarray*}
and
\begin{eqnarray*}
|\nabla w_{1}|^{\alpha}(F(D^{2}w_{1})+h(x)\cdot\nabla w_{1})+c(x)|w_{1}|^{\alpha}w_{1}-k(w_{1}+\delta)^{1+\alpha}= f_{\delta}(x,w_{0}).
\end{eqnarray*}
Using the comparison Theorem  \ref{Comparisontheorem} (with  $w_{0}=w_1$ on $\partial\Omega$)  we have that $w_{1}\geq w_{0}$.\\
Suppose that $w_{n}\geq w_{n-1}$ and let us show that $w_{n+1}\geq w_{n}$:\\
%
%\begin{eqnarray*}
%\Rightarrow\ |\nabla w_{0}|^{\alpha}(F(D^{2}w_{0})+h(x)\cdot\nabla w_{0})+c(x)|w_{0}|^{\alpha}w_{0}\geq 0\\
%\end{eqnarray*}
Since $f_{\delta}$ is decreasing, we have that
\begin{eqnarray*}
|\nabla w_{n}|^{\alpha}(F(D^{2}w_{n})&+&h(x)\cdot\nabla w_{n})+c(x)|w_{n}|^{\alpha}w_{n}-k(w_{n}+\delta)^{1+\alpha}\\
&=&f_{\delta}(x,w_{n-1})\geq f_{\delta}(x,w_{n}).
\end{eqnarray*}
This implies that 
\begin{eqnarray*}
|\nabla w_{n}|^{\alpha}(F(D^{2}w_{n})&+&h(x)\cdot\nabla w_{n})+c(x)|w_{n}|^{\alpha}w_{n}-k(w_{n}+\delta)^{1+\alpha}\\
&\geq& |\nabla w_{n+1}|^{\alpha}(F(D^{2}w_{n+1})
+h(x)\cdot\nabla w_{n+1})\\
&+&c(x)|w_{n+1}|^{\alpha}w_{n+1}-k(w_{n+1}+\delta)^{1+\alpha},
\end{eqnarray*}
and using the comparison Theorem \ref{Comparisontheorem} (with  $w_{n}=w_{n+1}$ on $\partial\Omega$) one gets   $w_{n+1}\geq w_{n}$. 

We have  shown that $\{w_{n}\}$ is non decreasing and since  $w_{0}>0$ in $\Omega$ one gets  $w_{n}> 0$  in $\Omega$ for all $n\geq 0$.  

Using the fact that $u^\star$ satisfies 
$|\nabla u^\star|^{\alpha}(F(D^{2} u^\star )+h(x)\cdot\nabla u^\star )+c(x)|u^\star |^{\alpha}u^\star-k(u^\star +\delta)^{1+\alpha}\leq  f_{\delta}(x,u^\star)\leq f_\delta ( x, w_n)$  we get at each step that $w_{n+1} \leq u^\star$, once we have assumed that $w_n \leq u^\star$.

Since the sequence  $\{w_{n}\}$ satisfies  the Lipschitz  estimates recalled  in Theorem \ref{lip} it converges uniformly to a function $Z_\delta$  which satisfies 
$$
|\nabla Z_\delta|^{\alpha}(F(D^{2}Z_\delta)+h(x)\cdot\nabla Z_\delta)+c(x)|Z_\delta|^{\alpha}Z_\delta+p(x)(Z_\delta+\delta)^{-\gamma}= 0.
$$
 Furthermore for any $\delta $ one has 
 $ u_\star \leq Z_\delta \leq u^\star$.

\subsubsection*{Step2: $\delta$ tends to $0$} 

Let $\delta$ and $Z_\delta$ defined by the first step. We note that since $u_\star \leq Z_\delta
\leq u^\star$, the term $p(x) (Z_\delta+ \delta)^{-\gamma}$ is uniformly locally  bounded  independently on $\delta$, and then $Z_\delta$ is uniformly locally Lipschitz.  It follows, using the uniform Lipschitz estimates in Theorem \ref{lip},  that one can extract from $Z_\delta$ a sequence which converges locally uniformly to some function $Z$,  such that $u_\star \leq Z
\leq u^\star$. Passing to the limit with Theorem \ref{Coseqlips}, one gets that, since $p( Z_\delta+ \delta)^{-\gamma}$ converges locally uniformly (for a subsequence) towards $pZ^{-\gamma}$, $Z$ is a solution of equation \eqref{eqfln} .
\subsection{Comparison principle and Uniqueness result}
We begin to prove some Lipschitz estimate between sub- and super-solutions of equation (\ref{eqfln})

  \begin{theorem}\label{lipsc}
     Suppose that $u$ is a positive,   bounded by above,  solution of 
     $$ \left\{\begin{array}{lc}
     | \nabla u |^\alpha ( F( D^2 u) + h(x) \cdot \nabla u) + p(x) u^{-\gamma} \geq  f &\ {\rm in}  \ \Omega \\
     u=0&\  {\rm on} \ \partial \Omega\end{array} \right.$$
     and $v$ is a  positive solution of 
      $$  \left\{\begin{array}{lc}| \nabla v |^\alpha ( F( D^2 v) + h(x) \cdot \nabla v) + p(x)  v^{-\gamma} \leq  g&\ {\rm in}  \ \Omega \\
     v=0 &\  {\rm on}\  \partial \Omega\end{array} \right.$$
      with $f$, $g$ and $h$  continuous and bounded, and $p>0$ is Holder continuous of exponent $\tau_p$. Then  : 
     \begin{itemize}
     \item If $\sup_{ \overline{\Omega}}  ( u-v) >0 $, there exists $c$  depending on $\Omega, |u|_\infty, |f|_\infty, |g|_\infty, |h|_\infty , |p|_\infty $ so that for all $(x,y)$ in $\overline{\Omega}^2$
     $$ u(x)-v(y) \leq \sup ( u-v) + C |x-y|$$
     
     \item If $\sup_{ \overline{\Omega}}  ( u-v)=0$, and if there exist $\tau_1\leq \inf (1,  {2 \over  1+\gamma}) $, and $C_1>0$  so that 
     $$u( x) \leq C_1 d( x, \partial \Omega)^{\tau_1},$$
    then,  there exists some constant $C$  depending on $\Omega, |u|_\infty, |f|_\infty, |g|_\infty, |h|_\infty , \tau_p$, and $C_1$ so that  for all $(x,y)$ in $\overline{\Omega}^2$
      $$  u(x)-v(y) \leq   C |x-y|^\tau, $$
      with $\tau = \inf ( \tau_1,  { 2+\alpha+\tau_p \over 1+\alpha+ \gamma})$. 
      \end{itemize}
      \end{theorem}
      
       \begin{remark} 
        In the sequel, we will apply this result with $f = -c u^{1+\alpha} $ and $g= -c v^{1+\alpha} ,  $
   mainly to prove the uniqueness. On the other hand, we get H\"older regularity of this solution,  by the second part of the Theorem above,  by recalling that from  the first  sections, one has an exponent $\tau_1$ which can be taken arbitrarily close to $1$ in the case $\gamma <1$ and  is equal to ${2\over 1+\gamma}$ if $\gamma >1$. 
     
       \end{remark}
     \textbf{Proof}
     
      Let $\omega$ be  defined on $\mathbb{R}^+$ by 
      $$ \omega( s) = s- {s^{1+\varepsilon} \over 2(1+\varepsilon)}$$
	      where $\varepsilon \in ]0,1[$. 
       Let us introduce  
       $$ \psi(x,y) = u(x)-v(y) -\sup ( u-v) -M \omega( |x-y|)$$
        where $M$ will be chosen large enough later. 
        
        It is clear that it is sufficient to prove that for $|x-y| < {1\over 2}$, 
        $$\psi(x,y) \leq 0$$
        Indeed, if $|x-y| > {1\over 2}$ and  if we assume that  ${M\over 2} > \sup u-\inf v$ the required result holds.  
        
         We argue by contradiction and suppose that $ \sup_{(x,y) \in \overline{ \Omega}^2} \psi( x, y) >0$. Then by the upper-semicontinuity of $\psi$,  it is achieved on some pair $(\bar x,\bar y) \in \overline{ \Omega}^2$.  In the following $\delta>0$ is a positive parameter, take 
         $$ M= {2(\sup u-\inf v)  \over \delta},$$ then from the definition of $\bar x$ and $\bar y$, 
         $| \bar x-\bar y | \leq \delta$.  So saying that $M$ is large is equivalent to say that $\delta$ is small. 
         
          We first remark that neither $\bar x$ nor $\bar y$ belongs to the boundary. Indeed, $\bar x$ cannot  be on the boundary by the positivity of $v$, and  if $\bar y \in \partial \Omega$ then we would have 
          \begin{equation}\label{eq1} u( \bar x) \geq \sup ( u-v) + {M \over 2} |\bar x-\bar y| .
          \end{equation}
         This is contradicted  by  the continuity of $u$, since  there exists $\delta_1$ so that for $d(x, \partial \Omega)
< \delta_1$ one has 
$$u(x) \leq    {\sup ( u-v) \over 2}$$
 while if $d( x, \partial \Omega) > \delta$
 $${ M\over 2 } d( \bar x,  \partial \Omega)+ \sup ( u-v)  \geq { M \delta \over 2} > \sup u $$
  which also contradicts (\ref{eq1}). 
  
   We have obtained that $(\bar x, \bar y) \in \Omega^2$. Furthermore $\bar x \neq \bar y$. 
        
           By Ishii's lemma (Lemma 9) in \cite{I} (see also \cite{BD3}) for all $\zeta >0$ there exist $X_\zeta $ and $Y_\zeta $ so that 
          
    $$(q, X_\zeta ) \in \overline{J}^{2,+} u( \bar x), \ ( q, -Y_\zeta  ) \in  \overline{J}^{2,-} v( \bar y)$$
    with       $q = M  \omega^\prime ( |  \bar x-\bar y|)$,      and      $$ \left( \begin{array}{cc} X_\zeta   & 0\\
     0& Y_\zeta   \end{array}\right) \leq M  \left( \begin{array}{cc} B & -B\\ -B& B\end{array} \right) +\zeta  M^2\left( \begin{array}{cc} B^2 & -B^2\\ -B^2& B^2\end{array} \right) $$
      with 
      $B = D^2( \omega ( |\cdot | ) ( \bar x-\bar y)$. One has $$B = \left( \omega^{\prime \prime }-{\omega^\prime \over  r}\right) {x\otimes x \over |x|^2}+ { \omega^\prime  \over r} I$$
             and then 
             $$B^2 =  \left(( \omega^{ \prime \prime} )^2-{ (\omega^\prime)^2(r) \over r^2}\right) {x\otimes x \over |x|^2} + { (\omega^\prime)^2(r) \over r^2} I.$$
             So taking    $\zeta = {1 \over  M(1+ 2 | \omega^{\prime \prime} (r) | + { \omega^\prime (r) \over r})}$,               $B+ \zeta  M  B^2$ has the eigenvalues 
                $\omega^{ \prime \prime} + \zeta  M  (\omega^{ \prime \prime} )^2 \leq { \omega^{ \prime \prime} \over 2}$
                  and  ${ \omega^\prime (r) \over r}+ \zeta { (\omega^\prime)^2 \over r^2}\leq 2{ \omega^\prime \over r}$.
                With that choice of $\zeta,  $ dropping the index $\zeta $ for $X_\zeta $ and $Y_\zeta $,  one has       $ X+Y \leq 0$
    and  for any $x$ 
        $$  ^t( x, -x) \left( \begin{array}{cc} X & 0\\
     0& Y \end{array}\right)  \left(\begin{array}{c}
     x\\
     -x
     \end{array}\right)\leq 4M ^tx (B +\varepsilon B^2) x.$$
  In particular there exist at least one eigenvalue which is  less or equal to $4  { \omega^{ \prime \prime} \over 2}  = 2 \omega^{ \prime \prime} $. 
  Using $X+ Y \leq 0$, one has 
  $tr(X+ Y) \leq 2M  \omega^{ \prime \prime}  $.   
             
   We have  then 
      $$tr( X+ Y) \leq -CM | \bar x-\bar y |^{ \varepsilon-1}$$
      while always by Ishii's lemma 

      $$| X| + |Y| \leq C  M  |B|_\infty \leq C M | \bar x-\bar y |. $$
      On the other hand       $$ {M\over 2} \leq |q| \leq M.$$

      Then, using in the following lines $F(X)-F( -Y) \leq a tr(X+Y)$, 
           \begin{eqnarray*}
      -|p|_\infty ( \sup ( u-v))^{-\gamma} + f( \bar x) &\leq &  -p(\bar x) ( u( \bar x) ^{ -\gamma} + f( \bar x)\\
       &\leq&  |q|^\alpha F( X) + h( \bar x) \cdot q | q |^\alpha \\
       &\leq & |q|^\alpha F( -Y) +  h( \bar y) \cdot q | q |^\alpha + 2| h|_\infty \ M^{1+\alpha}   \\
       &+ &M^{1+\alpha} ({1 \over 2})^ {| \alpha|}a tr( X+Y) \\
       &\leq & -p( \bar y) v( \bar y) ^{ -\gamma} +  g( \bar y) +2| h|_\infty \ M^{1+\alpha}  \\
       &-& C M^{1+\alpha} | \bar x -\bar y |^{\varepsilon -1}\\
       &\leq & 
      2 |h|_\infty\ M^{1+\alpha}  - C M^{1+\alpha} | \bar x -\bar y |^{\varepsilon -1}+ g( \bar y)
       \end{eqnarray*}
         from this we get  for some  positive constant  $C$, and for $\delta$ small  ( so that $| h|_\infty \ \delta^{1-\varepsilon } < { C \over 2}$), 
        $${C\over 2}  M^{1+\alpha} | \bar x-\bar y |^{ \varepsilon-1} \leq |p|_\infty  ( \sup ( u-v))^{-\gamma}+ |f|_\infty + |g|_\infty+ 2 |h|_\infty M^{1+\alpha}   , $$
        clearly a contradiction as soon as $M$ is large enough, since $\varepsilon <1$. We have then obtained that $\psi(x,y) \leq 0$ for all $x, y\in \overline{ \Omega}^2$.
        
         We now do the case where $\sup ( u-v) = 0$. We take the function 
         $$ \psi(x,y) = u(x)-v(y)-M | x-y|^\tau$$
          where $M =  {2(\sup u -\inf v)   \over \delta^\tau}$ , this  forces $\bar x$ and $\bar y$ to satisfy   
          $ | \bar x-\bar y | \leq \delta$.  We take $\tau = \inf ( \tau_1,  { 2+\alpha+\tau_p \over 1+\alpha+ \gamma})$. 
           We argue by contradiction and suppose that the supremum of $\psi$ is $>0$. Then it is achieved on some pair $( \bar x, \bar y) \in \overline{ \Omega}^2$. By taking $M$ larger than $C$ where $C$ is so that 
           $$u( \bar x) \leq C d(\bar  x, \partial \Omega)^{\tau_1}$$
            one obtains that $\bar y$ cannot belong to the boundary. 
             On the other hand, the positivity of $v$ implies that $\bar x$ cannot be on the boundary. 
             
                 Note for further purpose  that 
      $$ -p( \bar x) u( \bar x)^{-\gamma} \geq -p( \bar y)u( \bar x)^{-\gamma}  -C_p| \bar x-\bar y|^{\tau_p} ( M |\bar x-\bar y|^\tau)^{-\gamma}  \geq -p( \bar y) v( \bar y)^{-\gamma} - C M^{-\gamma} | \bar x-\bar y |^{ \tau_p-\tau \gamma}.$$
             
              We have then by Ishii's lemma,   for all $\zeta $ the existence of $X_\zeta  $ and $Y_\zeta  $ in ${ \cal S}$ so that with 
              $ q = M \tau |\bar x-\bar y |^{ \tau-1}$
              $$(q, X_\zeta ) \in \overline{J}^{2,+} u( \bar x), \ ( q, -Y_\zeta ) \in  \overline{J}^{2,-} v( \bar y)$$
         $$ \left( \begin{array}{cc} X_\zeta  & 0\\
     0& Y_\zeta   \end{array}\right) \leq M \left( \begin{array}{cc} B+\zeta  M  B^2 & -B-\zeta   M B^2\\ -B+\zeta   M B^2& B+\zeta  M B^2\end{array} \right) $$
      with 
      $B = D^2( |\cdot |^\tau) ( \bar x-\bar y)$.  By  the choice of $\zeta  $ as  the first part of the proof, and dropping the index $\zeta $,  one has  $X+ Y \leq 0$, and $tr(X+Y) \leq -C M | \bar x-\bar y |^{ \tau-2}$. \\
                  We have then by using the fact that $u$ and $v$ are respectively sub-and super-solutions
         \begin{eqnarray*}
          -p( \bar y) v( \bar y)^{-\gamma} &-& C M^{-\gamma} | \bar x-\bar y |^{ \tau_p-\tau \gamma} + f( \bar x) \\
          &\leq & |q|^\alpha F( X) + M^{1+\alpha} (h( \bar x) \cdot  \bar x-\bar y) | \bar x-\bar y|^{ (\tau-1)\alpha-1} \\
          &\leq & |q|^\alpha F( -Y) +  M^{1+\alpha} |h|_\infty | \bar x-\bar y|^{ (\tau-1)\alpha}  -C M^{1+\alpha } | \bar x-\bar y |^{(\tau-1)\alpha + \tau-2} \\
          &\leq &-p(\bar y) ( v( \bar y)^{-\gamma} + g( \bar y)+ M^{1+\alpha} |h|_\infty | \bar x-\bar y|^{ (\tau-1)\alpha}  -C M^{1+\alpha } | \bar x-\bar y |^{(\tau-1)\alpha + \tau-2}.
          \end{eqnarray*}
           From this one derives that 
           for some constants 
           $$ C M^{1+\alpha}  | \bar x-\bar y |^{(\tau-1)\alpha + \tau-2} \leq   C M^{1+\alpha} |h|_\infty | \bar x-\bar y|^{ (\tau-1)\alpha}  +  C M^{-\gamma} | \bar x-\bar y |^{ \tau_p-\tau \gamma}+ |f|_\infty+ |g|_\infty$$
            which is a contradiction  as soon as $\delta $ is small enough, by the assumption on $\tau$. We have obtained that 
             for all $x, y$ in $\overline{ \Omega}
             $
              $$u(x)-v(y) \leq C | x- y |^\tau.$$
$\hfill\square$
\ \\

        \begin{corollary}
         The  solutions constructed in the proof of the previous section are  H\"older continuous up to the boundary,   with an  exponent $\tau$ arbitrary close to $1$
         when $\gamma <1$,   and $\tau\leq  { \tau_p + \alpha+2 \over 1+\alpha+\gamma}$, when $\gamma >1$. 
        They are in both cases  Lipschitz continuous inside $\Omega$.
        \end{corollary}
        
         \textbf{Proof}
          The Lipschitz  interior continuity is immediate by using the results of \cite{BD2}, remarking that  on a compact set of $\Omega$, by the strong maximum principle and the continuity of $u$,  $pu^{-\gamma}$ is 
         bounded. 
         $\hfill\square$
\ \\
Using Theorem \ref{lipsc} we have the following comparison result between sub- and super solutions.
       \begin{theorem} %\label{compa1}
     
        Suppose that $u, v>0$ are respectively sub- and super-solutions of 
        $$ | \nabla u |^\alpha (F( D^2 u)+ h(x)\cdot \nabla u)   + c(x) u^{1+\alpha} + p(x) u ^{-\gamma} = 0 \mbox{  in }\Omega$$
         and are zero on the boundary. Then 
         $u\leq v$ in $\Omega$. 
                \end{theorem}
\textbf{Proof.}
       Note that for $y<\varepsilon := \left( { \gamma\min p \over (1+\alpha) |c|_\infty}\right)^{1 \over 1+\alpha+ \gamma}$
        the function 
        $y\mapsto c(x) y^{1+\alpha} + p(x) y ^{-\gamma}$
         is decreasing.          By upper-semicontinuity, since $u$ is zero on the boundary,  for all $\varepsilon $ there exists $\delta$ so that for $d< \delta$ one has 
         $ u(x) \leq \varepsilon$ . 
       
          Suppose by contradiction that $u> v$ somewhere. We suppose first that the supremum of $u-v$ is achieved inside $d(x, \partial \Omega ) \leq  \delta$. Let $\bar x$ be  some point in this set where the supremum is achieved. Then one has 
          $0 < v(\bar x) < u(\bar x) \leq \varepsilon$. 
           In particular 
        $c(\bar x) u(\bar x) ^{1+\alpha} + p(\bar x) (u(\bar x)) ^{-\gamma}< c(\bar x) v(\bar x) ^{1+\alpha} + p(\bar x) (v(\bar x))^{-\gamma}$. 
          
           Using the usual doubling of variables, say defining  in the case $\alpha >0$
            $$ \psi_j(x,y) = u(x)-v(y)-{j\over 2} |x_j-y_j|^2$$
            (the case $\alpha <0$ requires the changes provided  at the end of the proof), 
             there exist  $x_j$ and $y_j$ in $\{x \in \Omega :\ d(x, \partial \Omega)  \leq  \delta\}$ and $(X_j ,Y_j)$ in ${\mathcal S}^2$ so that 
             $$(j(x_j-y_j), X_j) \in \overline{J^{2,+} }u( x_j),\ (j(x_j-y_j), -Y_j) \in \overline{J^{2,-} }v( y_j)$$
              and  
             $$\left( \begin{array}{cc}
              X_j&0\\
              0& Y_j\end{array} \right)\leq 2j \left( \begin{array}{cc}
              I&-I\\-I& I\end{array}\right).$$
              Furthermore by the boundary conditions, neither $x_j$, nor $y_j$ belong to $\partial \Omega$. 
              
   From the Lipshitz estimates between sub and super-solutions in Theorem \ref{lipsc}, one has 
                $$ j |x_j-y_j|^2 + \sup ( u-v) \leq  u(x_j)-v(y_j) \leq \sup ( u-v) + |x_j-y_j|.$$
                 From what we derive that 
                 $j|x_j-y_j|$ is bounded.  In particular, using $h$ continuous, one has 
                 $ |h(x_j)-h(y_j)| ( j |x_j-y_j|)^{1+\alpha}  = o(1)$.
                 
              Using  $X_j + Y_j \leq 0$, and, since $F$  satisfies \eqref{defF}, one can write 
             \begin{eqnarray*}
  -c(x_j) u(x_j) ^{1+\alpha} - p(x_j) {u( x_j)}^{-\gamma}  & \leq& | j( x_j-y_j)|^\alpha (F( X_j)+ h(x_j) \cdot j( x_j-y_j)) \\
  &\leq &| j( x_j-y_j)|^\alpha( F( -Y_j) + h(y_j) \cdot j( x_j-y_j))+o(1)\\
  &\leq &   -c(y_j) v(y_j) ^{1+\alpha} - p(y_j) (v( y_j)) ^{-\gamma}+o(1),  
  \end{eqnarray*}
   and  using the continuity of $c$, $h$  and $p$ and passing to the limit when $j$ goes to infinity, one gets a contradiction.
   
    Suppose now that $\bar x$ is not  in $\{x \in \Omega :\ d(x, \partial \Omega)  \leq  \delta\}$.  Then there exists $0< \kappa < M$ so that  for all $x\in \Omega$  such that $d(x) < \delta$, 
    $ u(x) \leq v(x) + \sup ( u-v) -\kappa:= v(x) + M-\kappa$. We then use the change of function $U = \log u $ and $V = \log (v+ M-\kappa)$ in the set 
     $\Omega_\delta = \{ x \in \Omega, d(x, \partial \Omega) > \delta\}$
      One has 
      $U \leq V$ on the boundary, and $U$ and $V$ are respectively sub-and super-solutions of the equations 
      
      $$| \nabla U |^\alpha (F( D^2 U+ \nabla U \otimes \nabla U)+ h(x) \cdot \nabla U) + c(x) + e^{-U(1+\alpha+\gamma)}  p(x)   \geq 0$$
                    and 
                     $$| \nabla V |^\alpha (F( D^2 V+ \nabla V \otimes \nabla V) + h(x) \cdot  \nabla V)+ c(x) + e^{-V(1+\alpha+\gamma) } p(x)     \leq 0.$$
                                          Using the comparison principle for these type of equations, (see for example \cite{BDL2}, Theorem 5.1),  remarking that 
                      $y\mapsto  e^{-y(1+\alpha+\gamma)} p(x) $ is decreasing one gets that 
                      $U \leq V$ everywhere in $\Omega_\delta$. This implies that 
                      $ u \leq v+ M-\kappa$ everywhere in $\Omega$, a contradiction with the definition of the supremum.

                                              We  breafly give  the changes to bring in the case $\alpha <0$. 
                         In that case in the first step, the function $\psi_j$ must be replaced by 
                         $$ \psi_j(x,y) = u(x)-v(y)-{j\over q} |x_j-y_j|^q$$ where $q > {\alpha +2 \over 1+\alpha}$.                           We next follow the lines in the comparison Theorem 3.3 in  \cite{BDL1},  (A key point  consists in observing that $x_j \neq y_j$). 
                      
                      $\hfill\square$

                       \begin{corollary}
                       
                        There is uniqueness of solution for the equation (\ref{eqfln}).                         \end{corollary}
                         
%  \textbf{Proof.}
%                           It is sufficient to apply the previous comparison principle  to $u>0$ and $v>0$ two solutions, and exchange them.  $\hfill\square$

\section{Regularity of the unique viscosity solution of (\ref{eqfln}).} 

\subsection{Interior regularity of the viscosity solution}\ \\
Interior regularity of the unique viscosity solution of the problem (\ref{eqfln}) is easily obtained by the results about  regularity of  viscosity solutions of the following equation
\begin{eqnarray}\label{eqflngen}
|\nabla u|^{\alpha}F(D^{2}u)+h(x)\cdot \nabla u  |\nabla u|^{\alpha}=f
\end{eqnarray}
 which are recalled in the introduction, ( \cite{IS}, \cite{BDL1}, \cite{BDCocv}). 
Indeed,   
let $u$ be a solution of equation \eqref{eqfln},  
let  $f=-c(x)u^{1+\alpha}-p(x)u^{-\gamma}$, and  let us consider the  equation 
 $$ 
  | \nabla v |^\alpha  ( F( D^2 v)  + h(x) \cdot \nabla v) = -c(x)u^{1+\alpha}-p(x)u^{-\gamma}\  {\rm in} \ \Omega. 
  $$

It is immediate to check that  $u $ is a viscosity solution of this equation, hence,  using the classical regularity results recalled below, $u\in {\mathcal C}^{1, \beta}( \Omega)$.
     
 \subsection{"Regularity" up to the boundary}  \ \\

 We begin to remark that when $\gamma >1$, even the Lipschitz regularity up to the boundary does not hold. Next we will see some cases in which we can ensure the ${\mathcal C}^1$ regularity for $\gamma <1$. 
We note that if $x_{0}\ \in\ \partial\Omega$,  $\overrightarrow{n}$ denotes the inner normal to $\partial\Omega$ at $x_{0}$,  and $\phi_1$ is some  eigenfunction for $\lambda_1^{ \beta_{c, \alpha, \gamma}}$, then  by Hopf boundary principle 
\begin{equation*}
\lim_{s\rightarrow 0^{+}}{\dfrac{\phi_{1}(x_{0}+s\overrightarrow{n})}{s}}=\lim_{s\rightarrow 0^{+}}{\dfrac{\phi_{1}(x_{0}+s\overrightarrow{n})-\phi_{1}(x_{0})}{s}}=\nabla \phi_{1}(x_{0})\cdot\overrightarrow{n}.
\end{equation*}
If $\gamma>1$, let us recall that for some convenient $b_1>0$ and   for  $t=\dfrac{2+\alpha}{1+\alpha+\gamma}<1$, $ b_{1}\phi_{1}(x)^{t}$ is a sub-solution, and then by the comparison principle,  $u(x)\geq\ b_{1}\phi_{1}(x)^{t}$. It follows that, for $s>0$,
\begin{equation*}
\dfrac{u(x_{0}+s\overrightarrow{n})-u(x_{0})}{s}\geq \ b_{1}\phi_{1}(x_{0}+s\overrightarrow{n})^{t-1}\dfrac{\phi_{1}(x_{0}+s\overrightarrow{n})}{s}.
\end{equation*}
Therefore since $t<1$, $\phi_1 = 0$ on the boundary, using the existence of a positive constant $C$ such that $\dfrac{\phi_{1}(x_{0}+s\overrightarrow{n})}{s}> C $  given by  Hopf principle, then
\begin{equation*}
\lim_{s\rightarrow 0^{+}}{\dfrac{u(x_{0}+s\overrightarrow{n})-u(x_{0})}{s}}=+\infty,
\end{equation*}
so $u$  cannot  be Lipchitz continuous  on  $\overline{\Omega}$.
\subsection{Regularity  up to the boundary when  $N=1$ }\ \\

Suppose that $N = 1$,  $h = c = 0$  and $ p\equiv 1$. We prove below the ${\mathcal C}^1$ regularity up to the boundary in the case $\gamma <1$,  while in the case $\gamma>1$ the solution cannot be Lipschitz continuous up to the boundary.  We can suppose without loss of generality that 
 $ \Omega = ]0,1[$. Let  us  consider the equation 
 $$ |u^\prime |^\alpha u^{\prime \prime} + u^{-\gamma} = 0,$$
  $$u(0) = u(1)=0,$$
   multiplying by $u^\prime$ and integrating one gets
   
   $${|u^\prime |^{ 2+\alpha} \over 2+\alpha} + { u^{1-\gamma} \over 1-\gamma} = C,$$
   where $C$ is a constant, and then when $\gamma >1$, $\lim_{ x\rightarrow (0,1)}  {|u^\prime |^{ 2+\alpha} \over 2+\alpha} =  +\infty$. 
   
For the special case where $\gamma=1$, multiplying by $u^\prime$ and integrating one gets the equation
     
   $${|u^\prime |^{ 2+\alpha} \over 2+\alpha} + \log{u} = C,$$
   Note that the equation is invariant by the change $x \mapsto 1-x$, $u$ is concave and  then ${1\over 2}$ is a maximum point for $u$,  so   $u^\prime({1\over 2} )=0$.    
Then, for some positive constant $C$ defined by  $ u({1\over 2} ) =  e^C $, 
         \begin{equation}\label{uprime} u^\prime (x) =\left\{\begin{array}{lc} \left((2+\alpha)(C-\log{u})\right)^{1\over 2+\alpha} & {\rm if} \ x<{1\over 2}\\
       -\left((2+\alpha) (C-\log{u})\right)^{1\over 2+\alpha} & {\rm if} \ x>{1\over 2}.
       \end{array}\right.
       \end{equation}
Consequently,  $$ \lim_{x\rightarrow 0}  u^\prime  = \infty.$$\\

        If  $\gamma  < 1$ the solutions are   given ,        for some positive constant $C$ defined by  $ u({1\over 2} ) =  (C  (1-\gamma))^{1\over 1-\gamma}$
         \begin{equation}\label{uprime} u^\prime (x) =\left\{\begin{array}{lc} \left((2+\alpha)(C-{u^{1-\gamma}\over 1-\gamma})\right)^{1\over 2+\alpha} & {\rm if} \ x<{1\over 2}\\
       - \left((2+\alpha)(C-{u^{1-\gamma}\over 1-\gamma})\right)^{1\over 2+\alpha} & {\rm if} \ x>{1\over 2}.
       \end{array}\right.
       \end{equation} 
             
            From equation \eqref{uprime}, $u^\prime$ is continuous up to the boundary. 
    \\
\subsection{Existence of radial solution when $N>1$  and regularity up to the boundary when $\gamma <1$ and $\Omega$ is a ball. }\ \\

In all this sub-section we  still suppose that $ h = c = 0$, $p\equiv 1$  and that $N \geq 2$. \\
 We begin to prove the existence of a radial solution in the particular case  where  $F = tr $. We will do the general case later.
   We begin to construct in a neighborhood of $0$ a solution by using  a fixed point argument. First of all  suppose $v(0)= 1, v^\prime (0) = 0$ (necessary in the radial case). Note that we follow the method employed in particular in \cite{DG}.  
   Let us consider the map $v \mapsto T(v)$ where 
   $$ T( v)(r)  = 1 -\int_0^r\left({1 \over s^{ (1+\alpha) ( N-1)}} \int_0^s  \lambda^{ (1+\alpha) ( N-1)} (1+\alpha) v^{-\gamma} ( \lambda) d\lambda \right)^{1\over 1+\alpha} ds.$$
     We prove that for $r_o$ small enough, $T$ possesses a fixed point defined in $[0, r_o[$.  
    Let us  define 
$$ r_o = \left({{((N-1)(1+\alpha)+1)^{1\over 1+\alpha}}(2+\alpha)\over 2^{1+ { | \alpha|+1 \over 1+\alpha} \gamma}(\max(\gamma,1+\alpha)) ^{1\over 1+\alpha}(1+\alpha)}  \right)^{1+\alpha \over 2+\alpha} ,  $$   and let us consider  the ball 
     
      $$B= \left\lbrace v \in {\mathcal C} ( B(0,r_o)), |v(x)-1|_\infty  < {1\over 2} \right\rbrace.$$
        Then for $r< r_o$, $T$ maps $B$ into itself. Indeed :  
   
\begin{eqnarray*}   
 |T(v)-1|_{\infty}&\leq& \int_0^r\left({1 \over s^{ (1+\alpha) ( N-1)}} \int_0^s  \lambda^{ (1+\alpha) ( N-1)} (1+\alpha) 2^\gamma  d\lambda \right)^{1\over  1+\alpha}ds\\
 &=& r^{2+\alpha \over 1+\alpha} {(2^\gamma(1+\alpha)) ^{1\over 1+\alpha}(1+\alpha) \over {((N-1)(1+\alpha)+1)^{1\over 1+\alpha}}(2+\alpha)}\\
&\leq&{1\over 2}.
\end{eqnarray*} 
         In order to check  that     $T$ is a contracting mapping on $B$,  
      we denote  for $v, w$ in $B$ 
 $$X(s) = {1\over s^{ ( 1+\alpha)(N-1)} }\int_0^s \lambda^{ ( 1+\alpha)(N-1)}  (1+\alpha) v^{-\gamma} d\lambda$$
  and 
   $$Y(s) = {1\over s^{ ( 1+\alpha)(N-1)}} \int_0^s \lambda^{ ( 1+\alpha)(N-1)}  (1+\alpha) w^{-\gamma} d\lambda.$$
    Note that 
    \begin{eqnarray*}
    \rm{i)}&&  | v^{-\gamma} - w^{-\gamma} |  \leq 2^{\gamma+1}\gamma  |v-w|,\\
    \rm{ii)}&& { {s 2^{-\gamma}(1+\alpha) \over  (1+\alpha) (  N-1)+1}\leq X(s)\leq {s 2^{\gamma}( 1+\alpha)  \over  (1+\alpha) (  N-1)+1}}\\
    && \ \ \ \  \rm{\ and\ }   {s 2^{-\gamma}1+\alpha) \over  (1+\alpha) (  N-1)+1}\leq Y(s) \leq {s 2^{\gamma}( 1+\alpha)  \over  (1+\alpha) (  N-1)+1},\\
    \rm{iii)}&&  |X(s)-Y(s)| \leq {1+\alpha \over (1+\alpha)(N-1)+1} |v^{-\gamma}-w^{-\gamma} | s.
    \end{eqnarray*}
    
      Then using  i), ii), iii) and the mean value's Theorem, for some $\theta \in ]0,1[$ 
      \begin{eqnarray*}
       | X(s)^{1\over 1+\alpha}-Y(s)^{1\over 1+\alpha}| &\leq& {1 \over 1+\alpha} | X(s)-Y(s)| |  X(s)+ \theta ( Y(s)-X(s))| ^{-\alpha\over 1+\alpha}\\
            &\leq &2^{1+ { | \alpha|+1 \over 1+\alpha} \gamma}  \gamma^{1\over 1+\alpha}\left( {1+\alpha \over (1+\alpha)(N-1)+1} \right)  ^{1 \over 1+\alpha} s^{1 \over 1+\alpha} |v-w|.
        \end{eqnarray*}
         As a consequence 
         \begin{eqnarray*}
          |T(v)-T(w)| &  \leq  & \int_0^r |X(s)^{1\over 1+\alpha} -Y(s) ^{1\over 1+\alpha} | ds \\
          &\leq & 2^{1+ { 1+| \alpha| \over 1+\alpha} \gamma} \gamma^{1\over 1+\alpha}   \left( {1+\alpha \over (1+\alpha)(N-1)+1} \right)  ^{1 \over 1+\alpha}\left({1+\alpha \over 2+\alpha}\right)  r^{2+\alpha \over 1+ \alpha} |v-w|,
          \end{eqnarray*} 
           and then $T$ is a contraction  mapping.

            This gives the local existence and uniqueness of a fixed point, denoted  $u$,   around $0$.  We can suppose, up to  replace $r_o$ by some smaller number, that $u^\prime  <0$  and $u>0$ in the whole interval $]0, r_o[$. Now if $ r_1 >0$  is so that $u^\prime ( r_1) <0$,   Cauchy Lipschitz Theorem gives the  local existence and uniqueness of a solution.  For that it is sufficient to consider the ordinary differential equation :
             $$\left( \begin{array}{c} 
        v^\prime\\ w^\prime \end{array}\right) := \varphi ( v, w) = \left( \begin{array}{c}
         |w|^{-{\alpha\over \alpha +1}} w \\
        ( - v^{-\gamma}-{N-1\over r}w )
          \end{array}\right),$$
          with $v(r_1) = u(r_1), \ w(r_1) = u^\prime ( r_1)$, 
           and to observe that $\varphi$   is a  Lipschitz function  of $(v, w)$ as long as neither $w$, nor $v$ takes the value $0$. 
           We denote by $u$ the  fixed point for $T$ in $[0, r_1]$ extended by the unique local solution of the previous  ODE,   as long as it is defined. 

           \noindent Note that    $u$ satisfies  on $[r_1, r[$
            $$ u(r) = u(r_1) - \int_{r_1}^r\left( {1 \over s^{ (1+\alpha) ( N-1)}} \int_0^s  \lambda^{ (1+\alpha) ( N-1)} (1+\alpha) u^{-\gamma} ( \lambda) d\lambda \right)^{1\over 1+\alpha} ds.$$
                         Then as long as $u >0$ one has  
             $$u(r) = 1- \int_0^r \left({1 \over s^{ (1+\alpha) ( N-1)}} \int_0^s  \lambda^{ (1+\alpha) ( N-1)} (1+\alpha) u^{-\gamma} ( \lambda) d\lambda \right)^{1\over 1+\alpha} ds.$$
              Since $u $ has values in $ ]0,1[$, and $u$ is not identically  equal to $1$
              $$\int_0^r \left({1 \over s^{ (1+\alpha) ( N-1)}} \int_0^s  \lambda^{ (1+\alpha) ( N-1)} (1+\alpha) u^{-\gamma} ( \lambda) d\lambda \right)^{1\over 1+\alpha} ds\geq    r^{2+\alpha  \over 1+\alpha} {(1+\alpha) ^{2+\alpha\over 1+\alpha} \over{ ((N-1)(1+\alpha)+1)^{1\over 1+\alpha}}(2+\alpha)} $$ and then  taking  $R$ so that 
              $ R^{2+\alpha \over 1+\alpha} {(1+\alpha) ^{2+\alpha\over 1+\alpha} \over {((N-1)(1+\alpha)+1)^{1\over 1+\alpha}} (2+\alpha)}> 1$, 
             one obtains  that there exists  $\bar r< R$,   so that  $u( \bar r) = 0$.

                      We  then  consider  $\tilde u$  defined as : 
              $$\tilde u ( r) = C u( \bar r \ r), \rm{\ with} \ C := \bar r^{-{ 2+\alpha\over  \gamma + \alpha +1}}.$$
                Then 
                $\tilde u $ solves the equation 
                $$ | \tilde u^\prime |^\alpha ( \tilde u^{ \prime \prime} + { N-1\over r} \tilde u^\prime ) +  \tilde u^{-\gamma} = 0$$
                 and $\tilde u( 1) = 0$. 
                 
                  The computations  above can easily be generalized to the case where $p$ is a radial function which satisfies the assumptions of the article.

We now observe that in the radial case, when $\gamma <1$ the solution $\tilde u$ above is ${\mathcal C}^1$.  Indeed, multiplying 
  $$|\tilde u^\prime |^\alpha \tilde u^{ \prime \prime} + {N-1\over r}  |\tilde u^\prime |^\alpha \tilde u^\prime  + \tilde u^{-\gamma} = 0$$
   by $ \tilde u^\prime r^{ (N-1)(2+\alpha)}$ and integrating, one has 
    $$ { d \over dr} \left( \frac{|\tilde u^\prime |^{ 2+\alpha} }{2+\alpha}r^{ (N-1) (2+\alpha)}+ {\tilde u^{1-\gamma}  \over 1-\gamma} r^{ (N-1) (2+\alpha)}\right) = (N-1)(2+\alpha) r^{  (N-1) (2+\alpha)-1} {\tilde u^{1-\gamma}\over 1-\gamma} ,$$
     and then, integrating between ${1\over 2}$ and $r$ one gets that  
      $$ \frac{|\tilde u^\prime |^{ 2+\alpha} }{2+\alpha} r^{ (N-1) (2+\alpha)}+ {\tilde u^{1-\gamma}  \over 1-\gamma} r^{ (N-1) (2+\alpha)}-C={1\over 1-\gamma}  \int_{1\over 2}^r (N-1)(2+\alpha) s^{  (N-1) (2+\alpha)-1} \tilde u^{1-\gamma}(s) ds,$$
   where 
   $C = \left( \frac{|\tilde u^\prime |^{ 2+\alpha} }{2+\alpha}r^{ (N-1) (2+\alpha)}+ {\tilde u^{1-\gamma}  \over 1-\gamma} r^{ (N-1) (2+\alpha)}\right)({1\over 2})$,    which proves that $ |\tilde u^\prime |^{ 2+\alpha} $ has  a finite limit when $r$ goes to $1$.
   
   We now do the general case .     We begin with  the case of  one of the Pucci's operator. We will deduce the general case by using the fact that the operator $F$ is sandwiched between  the two Pucci's operators, ${\mathcal M}^+_{a, A}$ and ${\mathcal M}^-_{a, A}$.

    Suppose that $F = {\mathcal M}^+_{a, A} $. We begin to prove   local existence and uniqueness of a solution near $0$. 
       For that aim we argue as in the case of the Laplacian , say we observe that  if $r_o$ is replaced by $r_o a^{1\over 2+\alpha}$ and $B$ is defined as in the  Laplacian case, we have a fixed point,  denoted  $u_o$, defined on $[0, r_o[$   for the operator  $T$  in $B$ : 
       $$ T(v) = 1-\int_0^r \left({1\over s^{(N-1)(1+\alpha)}} \int_0^s \lambda^{(N-1)(1+\alpha)} {(1+\alpha) \over a}v ^{-\gamma}d \lambda \right)^{1\over 1+\alpha}ds.$$Up to replace $r_o$ by some smaller number, one can assume that for $r < r_o$, ${u^{-\gamma}\over a}  + {(N-1) u^\prime \over r} < 0$. Let $r_1$ be so that $r _1 \in ]0, r_o[$.   We consider  for $r>r_o $, the ordinary differential equation 
        $$\left( \begin{array}{c} 
        v^\prime\\ w^\prime \end{array}\right) := \varphi ( v, w) = \left( \begin{array}{c}
         |w|^{-{\alpha\over \alpha +1}} w \\
          f_{a, A} ( - v^{-\gamma}-{N-1\over r}aw )
          \end{array}\right)$$
           where
           $$f_{ a, A} (x) = {x^+\over A} - {x^-\over a},$$
          with $v( r_o ) = u(r_o )\neq 0$, $w( r_o )= u^\prime (r_o )\neq 0$. The      function $\varphi$ is Lipschitz continuous as long as $v\neq 0$ and $w\neq 0$. Then Cauchy Lipschitz Theorem ensures  local  existence and uniqueness  of solution for $r> r_o $.  Let $r_1> r_o$ and $u_1$ defined on $[r_o, r_1[$  which solves this ordinary differential equation.  Let 
          $$ u= \left\{\begin{array}{cc}
           u_o& {\rm if} \ r < r_o\\
            u_1 & {\rm if} \ r \in [r_o, r_1[\end{array}\right.$$
               We observe that $u^\prime < 0$ as long as $u(r)>0$. Indeed, suppose that  $r_2> r_o$ is so that  $u^{\prime}(r)\leq 0$, $u(r)>0$  for $r < r_2$,  and $u^\prime (r_2) = 0$. One would have, since $ ( |u^\prime|^\alpha u^\prime ) $ is continuous,   $ ( |u^\prime|^\alpha u^\prime ) ^\prime(r_2)  \geq 0$, while by the equation this quantity is $<0$ on $r_2$, a contradiction.  As a consequence,       $u$ is a solution of the equation related to ${ \cal M}^+_{a, A}$ on $[0, R[$ where $R\leq \infty$ is so that  $u(r) >0$, $u^\prime(r)<0$   for $r< R$ and $\lim_{r\rightarrow R} \inf (u, u^\prime )(r) = 0$.
       As a conclusion we have obtained a solution $u$ on some intervall $[0, R[$, with $R$  so that $u(r) >0$, $u^\prime (r) <0$ for $r < R$,  and $\lim_{r\rightarrow R} u(r) = 0$.  In the following lines  we prove that $R < \infty$. 

     Since $u^{\prime}\neq 0 $ on $]0, R[$,  $u^{\prime\prime}$ is continuous and  then the sets $\{r>0,\ u^{\prime\prime}(r)<0\}$ and
$\{r>0,\ u^{\prime\prime}(r)>0\}$ are open. Then each of these sets is a countable union of intervals.     

      So there exist a numerable set of $r_i$ so that 
      on $]0, r_1[$,  and on $]r_{2i}, r_{2i+1}[ $,  $u^{\prime\prime } <0$, on $]r_{2i+1},  r_{2i+2} [$, $u^{\prime\prime } >0$.   Then on 
      $]r_{2i}, r_{2i+1}[ $,   $u$ satisfies 
      \begin{eqnarray*}
      -(-u^\prime )^{1+\alpha} ( s) &=&- (-u^\prime )^{1+\alpha} (r_{2i}) \left({r_{2i} \over s} \right)^{ (N-1)(1+\alpha)} \\
      &-&{1 \over s^{ (N-1)(1+\alpha)} } \int_{r_{2i}}^s {(1+\alpha) \over a}  u^{-\gamma } ( \lambda ) \lambda^{(N-1)(1+\alpha)} d\lambda\\
      &:=& h_{2i} (s)( \leq 0),
      \end{eqnarray*}
       and then 
    $$ u(r) = u(r_{2i}) - \int_{r_{2i}}^r (-h_{2i}(s))^{1\over 1+\alpha}  ds .$$
       On  $]r_{2i+1}, r_{2i+2}[$ we have the analogous formula with 
     
  \begin{eqnarray*}
&h_{2i+1} ( s)& = - (-u^\prime )^{1+\alpha}(r_{2i+1})  \left({r_{2i+1} \over s} \right)^{ (N-1)(1+\alpha)a\over A}\\
    &   &\ \ \ \ -{1 \over s^{(N-1)(1+\alpha)a\over A} } \int_{r_{2i+1}}^s  {(1+\alpha)  \over A}  u^{-\gamma } (\lambda) \lambda^{(N-1)(1+\alpha)a\over A} d\lambda.    
\end{eqnarray*}    
      We  observe next that there exists $\bar r$ so that 
      $u(\bar r )  =0$. For that aim we write  for $r\in ]r_{2k-1}, r_{2k}[$
      $$ u( r) = u(0)+ (u(r_{0})-u(0))+ \sum_0^{k-1}\left( u(r_{2i+1})-u(r_{2i})\right)+ u(r)-u(r_{2k-1}).$$
        Note that $u(r_{2i+1})-u(r_{2i})<0$, and  since $u$ is bounded  by $u(0)$,
        \begin{eqnarray*}
         u(r)-u(r_{2k-1} )& \leq& -(u(0))^{ -\gamma\over 1+\alpha} \int_{r_{2k-1}}^r \left(  {1 \over s^{ (N-1)(1+\alpha)a\over A} } \int_{r_{2i+1}}^s  {(1+\alpha)  \over A} \lambda^{(N-1)(1+\alpha)a\over A} d\lambda\right)^{1\over 1+\alpha} ds\\
         &\leq& - C(u(0))^{ -\gamma\over 1+\alpha} (r^{2+\alpha \over 1+\alpha}-r_{2k-1}^{2+\alpha \over 1+\alpha} ).
         \end{eqnarray*}
                  where $C$ is a positive constant.          An analogous formula holds for $r \in ]r_{2k}, r_{2k+1}[$. 
         
Suppose that $r_{2i} >1$, then  for all $k\geq i$,  on $[r_{2k}, r_{2k+1}[$, 
$-\lambda^{(N-1)(1+\alpha)} \leq -\lambda^{(N-1)(1+\alpha)a \over A}$ so in the previous  inequality we do "as if " the same formula holds for the case $u^{\prime \prime} \geq 0$ or $\leq 0$ as soon as $r_{2i}$ is  greater than $1$, and then for $r \in [r_{2i}, r_{2i+1}]$ 
$$ u(r) \leq u(r_{2i} )-\int_{r_{2i}}^r \left({1 \over s^{ (N-1)(1+\alpha)} } \int_{r_{2i}}^s {(1+\alpha) \over a}  u^{-\gamma } ( \lambda ) \lambda^{(N-1)(1+\alpha)a\over A} d\lambda\right)^{1\over 1+\alpha} ds.$$
 Finally one has 
 for any $r$, $r>1$ , denoting by $i_o$ the first $i$ so that $r_{2(i_o+1)} >1$, if it exists, 
 \begin{eqnarray*}
  u(r)& \leq& u(0) + (u(r_{0})-u(0))+\sum_{i, i\leq i_o-1} ( u(r_{2i+1})-u(r_{2i}) )+  u(r)-u(r_{2i_o})\\
& \leq& u(0) - \int_{r_{2i_o}} ^r \left({1 \over s^{ (N-1)(1+\alpha)} } \int_{r_{2i}}^s {(1+\alpha) \over a}  u^{-\gamma } ( \lambda ) \lambda^{(N-1)(1+\alpha)a\over A} d\lambda\right)^{1\over 1+\alpha}ds \\
&\leq & u(0) - Cu(0)^{-\gamma\over 1+\alpha} (  r^{2+\alpha\over 1+ \alpha} - r_{2i_o}^{2+\alpha\over 1+ \alpha} )\\
&\leq & u(0) + C u(0)^{-\gamma \over 1+\alpha}  -Cu(0)^{-\gamma\over 1+\alpha}  r^{2+\alpha\over 1+ \alpha}. 
\end{eqnarray*}
 And then for $r$ large enough this quantity becomes negative. 
  If for all $i, r_{2i}\leq 1$, since $r_{2i}$ is increasing let $l$ its limit, then one can write  for $r> l$
 $$ u(r) \leq u(0)  - Cu(0)^{-\gamma\over 1+\alpha} (  r^{2+\alpha\over 1+ \alpha} -l^{2+\alpha \over 1+ \alpha}) \leq u(0) - Cu(0)^{-\gamma\over 1+\alpha} (  r^{2+\alpha\over 1+ \alpha}-1)  $$ and the same  conclusion follows.  
  Let then $\bar r$  be so that 
  $ u( \bar r)=0$.             
             We end the proof as in the case where ${\mathcal M}^+_{a, A}$ is replaced by the Laplacian, and we have obtained a  radial solution of the equation related to ${ \cal M}^+_{a, A}$ in the ball $B(0,1)$. 
             
                            For the general case we observe that  by the previous computations, there exists $\bar u
              $ a  radial solution for 
              $$|\bar u^{\prime}|^{\alpha}{\mathcal M}^+_{a, A}( D^2 \bar u)= -\bar u^{-\gamma}, \bar u( 1) = 0. $$
      Then it  provides   a super-solution for the equation, while, by obvious changes in the analysis above,   there exists  $\underline{u}$ 
               a radial solution for 
               $$ |\underline u^{\prime}|^{\alpha}{\mathcal M}^-_{a, A}( D^2 \underline{u} ) = -{ \underline u} ^{ \gamma}, \underline{u}(1) = 0$$
                which provides 
                a sub-solution. By the comparison principle (  in the uniqueness part), $\underline{u} \leq \overline{u}$. Using  Perron's method adapted to the present context (see \cite{BD1}), we obtain the existence of a radial  solution of (\ref{eqfln})
                  which lies  between $\bar u$ and $\underline{u}$. \\
                  
 { \em \textbf{Acknowledgment}}  : The author wishes to thank the anonymous referee for his  remarks which permit to improve considerably this paper.

 {\em \textbf{Email adress} : cheikhou-oumar.ndaw@cyu.fr}

\begin{thebibliography}{99}






 \bibitem{ART} D. J. Araújo, G. Ricarte, E. Teixeira, {\em Geometric gradient estimates for solutions to degenerate elliptic equations}. Calc. Var. Partial Differential Equations 53 (2015), no. 3-4, 605-625. 
 
 
   \bibitem{AR}A. Attouchi, E. Ruosteenoja, {\em Remarks on regularity for $p$Laplacian type equations in non divergence form}, Journal of Diff. Equations (2018), {265} (5), 1922-1961.
   \bibitem{BNV} 
   H. Beresticki, L. Nirenberg, S.R.S. Varahdan, {\em The principal eigenvalue and maximum principle for second order elliptic operators in general domain}, Comm. Pure Appl. Math., 47 (1994), no. 1, 47-92.
  
   \bibitem{BD1} I. Birindelli, F. Demengel, \emph{First eigenvalue and Maximum principle for fully nonlinear singular
operators},   Advances in Differential equations,
(2006) {11} (1), 91--119. 


 \bibitem{BD2}I. Birindelli, F. Demengel, {\em  Fully nonlinear operators with Hamiltonian: Hölder regularity of the gradient}. NoDEA Nonlinear Differential Equations Appl. 23 (2016), no. 4, Art. 41, 17 pp.
 
 \bibitem{BD3}I. Birindelli, F. Demengel, {\em  Existence and regularity results for fully nonlinear operators on the model of the pseudo Pucci's operators}.  J. Elliptic Parabol. Equ. 2 (2016), no. 1-2, 171–187.

\bibitem{BDjde} I. Birindelli, F. Demengel, {\em Regularity and uniqueness of the first eigenfunction for singular fully nonlinear operators}. J. Differential Equations 249 (2010), no. 5, 1089–1110. 

  
   \bibitem{BDCocv}I. Birindelli, F. Demengel, {\em ${\mathcal C}^{1, \beta} $regularity for Dirichlet problems associated to Fully Non linear equations degenerate elliptic equations} , COCV, (2014) , 1009-1024. 
   
      \bibitem{BDL1}  I. Birindelli, F. Demengel, F. Leoni,  {\em Dirichlet problems for fully nonlinear equations with "subquadratic" Hamiltonians}. Contemporary research in elliptic PDEs and related topics, 107-127, Springer INdAM Ser., 33, Springer,  2019. 
  \bibitem{BDL2}  I. Birindelli, F. Demengel, F. Leoni, {\em Ergodic pairs for singular or degenerate fully nonlinear operators},  ESAIM Control Optim. Calc. Var. 25 (2019).
  
\bibitem{BG}  I. Birindelli, G. Galise {\em The Dirichlet problem for fully nonlinear degenerate elliptic equations with a singular nonlinearity.} Calc. Var. 58, 180 (2019).
  
   \bibitem{C}  Luis A. Caffarelli, {\em Interior a priori estimates for solutions of fully nonlinear equations} . Ann. of Math. (2) 130 (1989), no. 1, 189-213. 
   \bibitem{CC}
Luis A. Caffarelli, X. Cabre,  Fully nonlinear elliptic equations. American Mathematical Society Colloquium Publications, 43. American Mathematical Society, Providence, RI, 1995.
\bibitem{CRT} M.G. Crandall,  P.H. Rabinowitz, L. Tartar, {\em On a Dirichlet problem with singular nonlinearity.} Commun. Partial Differ. Equ. 2, 193-222 (1977).
  
  \bibitem{usr} M.G. Crandall, H. Ishii, P.L. Lions, {\em User's guide to
viscosity solutions of second order partial differential equations}, Bull.
Amer. Math. Soc. (N.S.) 27 (1992), no. 1, 1--67.
\bibitem{DG} F. Demengel, O. Goubet, {\em Existence of boundary blow up solutions for singular or degenerate fully nonlinear equations},  Commun. Pure Appl. Anal. 12 (2013), no. 2, 621-645.
  
\bibitem{FQS} P. Felmer, A. Quaas, B. Sirakov {\em Existence and regularity results for fully nonlinear equations with singularities.} Math. Ann. 354 (2012), no. 1, 377-400.
  
 \bibitem{I}{H. Ishii}, {\em Viscosity solutions of non-linear partial differential equations.} Sugaku Expositions vol 9,  135-152 (1996).
 
 \bibitem {IL} {H. Ishii, P.L. Lions}, {\it  Viscosity solutions of Fully-Nonlinear Second  Order Elliptic Partial Differential Equations}, {J. Differential Equations},  {83},  (1990), 26--78.

   \bibitem{IS} C. Imbert, L. Silvestre, {\em ${\mathcal C}^{1, \alpha}$ regularity of solutions of some degenerate fully non-linear elliptic equations},  advances in Mathematics, (233), (2013), pp. 196-206. 
   
   
    \bibitem{LMK}  Lazer, McKenna, {\em On a singular nonlinear elliptic boundary-value problem}, Proc. Amer. Math. Soc. 111 (1991), 721-730 
    
\bibitem{V}  A. Vitolo {\em Singular elliptic equations with directional diffusion.} Mathematics in Engineering 3(3), 1-16 (2021). 
    
 
    \end{thebibliography}
          \end{document}